\documentclass{amsart}

\usepackage[hypertex]{hyperref}
\usepackage{hyperref}
\usepackage[capitalise]{cleveref}

\usepackage{color}

\renewcommand{\tilde}{\widetilde}

\newcommand{\pf}[1]{{\langle\!\langle{#1}\rangle\!\rangle}} %formes de
\newcommand{\ba}{\overline{\rule{2.5mm}{0mm}\rule{0mm}{4pt}}} %canonical involution

\newcommand{\cf}{\mathcal{F}}

\newcommand{\C}{\mathcal{C}}
\newcommand{\cd}{\mathcal{D}}

\newcommand{\mz}{\mathbb{Z}}
\newcommand{\mq}{\mathbb{Q}}

\newcommand{\qform}[1]{{\langle{#1}\rangle}} % formes quadratiques
\newcommand{\pform}[1]{{\langle\!\langle{#1}\rangle\!\rangle}} %formes de Pfister

\DeclareMathOperator{\br}{Br}
\DeclareMathOperator{\disc}{disc}

\DeclareMathOperator{\hyp}{hyp}
\DeclareMathOperator{\End}{End}
\DeclareMathOperator{\ad}{ad}

\DeclareMathOperator{\cores}{cor}
\DeclareMathOperator{\res}{res}
\DeclareMathOperator{\id}{id}
\DeclareMathOperator{\Nrd}{Nrd}
\DeclareMathOperator{\nrp}{Nrp}

\DeclareMathOperator{\Int}{Int}

\DeclareMathOperator{\Orth}{O}
\DeclareMathOperator{\GO}{GO}
\DeclareMathOperator{\SO}{O^+}
\DeclareMathOperator{\SOn}{O_n^+}
\DeclareMathOperator{\SL}{SL}
\DeclareMathOperator{\SU}{SU}
\DeclareMathOperator{\PGU}{PGU}
\DeclareMathOperator{\PGSp}{PGSp}
\DeclareMathOperator{\PGO}{PGO}
\DeclareMathOperator{\GU}{GU}
\DeclareMathOperator{\SB}{SB}
\DeclareMathOperator{\Spin}{Spin}
\DeclareMathOperator{\Sp}{Sp}
\DeclareMathOperator{\Sym}{Sym}
\DeclareMathOperator{\SSym}{SSym}

\newtheorem{lem}{Lemma}[section]
\newtheorem{prop}[lem]{Proposition}
\newtheorem{thm}[lem]{Theorem}
\newtheorem{cor}[lem]{Corollary}

\theoremstyle{remark}
\newtheorem{remark}[lem]{Remark}

\newtheorem{example}[lem]{Example}
\newtheorem{defi}[lem]{Definition}

\title[Degree $3$ relative invariant for unitary involutions]{Degree $3$ relative invariant for unitary involutions}   

\author{Demba Barry} 
\address{Facult\'e des Sciences et Techniques de Bamako, BP: E3206 Bamako, Mali  and Departement Wiskunde--Informatica, Universiteit Antwerpen, Belgium}
\email{Barry.Demba@gmail.com}
\thanks{The first author gratefully acknowledges support from the FWO Odysseus Programme (project Explicit Methods in Quadratic Form Theory).}

\author{Alexandre Masquelein}
\address{ICTEAM Institute, Box L4.05.01\\
Universit\'e catholique de Louvain\\
B-1348 Louvain-la-Neuve, Belgium}
\email{alexandremasquelein@gmail.com}
\thanks{Part of the results of this paper were contained in the second author's PhD dissertation, which was supervised by Jean-Pierre Tignol. }

\author{Anne Qu\'eguiner-Mathieu}
\address{
LAGA - CNRS (UMR 7539)\\
Universit\'e Sorbonne Paris Nord \\
F-93430 Villetaneuse, France}
\email{queguin@math.univ-paris13.fr}
\thanks{The first and third authors acknowledge support from the CIRM-Luminy and its Research in Pairs program.}

\subjclass[2010]{}

\date{\today}

\begin{document}

\begin{abstract}
Using the Rost invariant for non split simply connected groups, we define a relative degree $3$ cohomological invariant for pairs of orthogonal or unitary involutions having isomorphic Clifford or discriminant algebras. The main purpose of this paper is to study general properties of this invariant in the unitary case, that is for torsors under groups of outer type ${\sf A}$.  If the underlying algebra is split, it can be reinterpreted in terms of the Arason invariant of quadratic forms, using the trace form of a hermitian form. 
When the algebra with unitary involution has a symplectic or orthogonal descent, or a symplectic or orthogonal quadratic extension, we provide comparison theorems between the corresponding invariants of unitary and orthogonal or symplectic types. We also prove the relative invariant is classifying in degree $4$, at least up to conjugation by the non-trivial automorphism of the underlying quadratic extension.  
In general, choosing a particular base point, the relative invariant also produces absolute Arason invariants, under some additional condition on the underlying algebra. Notably, if the algebra has even co-index, so that it admits a hyperbolic involution, which is unique up to isomorphism, we get a so-called {\em hyperbolic} Arason invariant. Assuming in addition the algebra has degree $8$, we may also define a {\em decomposable} Arason invariant. It generally does not coincide with the hyperbolic Arason invariant, as the hyperbolic involution need not be totally decomposable.

\end{abstract}

\maketitle

\section{Introduction and notations}  

In quadratic form theory, the Arason invariant, introduced by Arason in~\cite{Arason}, is a cohomological invariant, with values in $H^3(F,\mu_2)$, associated to any $2m$-dimensional quadratic form with trivial discriminant and trivial Clifford invariant. It is closely related to the Rost invariant of the algebraic group $\Spin_{2m}$, see~\cite[p. 107--108  and Thm. 9.11]{GMS}. Using the Rost invariant for non split Spin groups, and for absolutely almost simple simply connected groups of other types, one may try to define analogues of the Arason invariant for the underlying algebraic objects, namely hermitian forms and involutions. This was initiated by Bayer and Parimala in~\cite{BP}. They defined a degree $3$ invariant for hermitian forms, which plays a role in their proof of the so-called Hasse Principle Conjecture II. 

For involutions of the first kind, degree $3$ cohomological invariants were investigated in several papers, such as~\cite{BMT} and ~\cite{GPT} for symplectic involutions, as well as~\cite{Ber}, ~\cite{G}, ~\cite{QT-deg8} and~\cite{QT-Arason} for orthogonal involutions; see also~\cite{Tignol:hyderabad} for a survey on these questions. Following a suggestion of Tignol~\cite[\S4.3]{Tignol:hyderabad}, the present paper investigates the case of unitary involutions, which correspond to groups of outer type ${\sf A}$. 

The first two sections are devoted to definitions and formal properties of analogues of the Arason invariant for unitary involutions. The results there do not depart much from what is known in the orthogonal case. In particular, the Arason invariant for quadratic forms may be used to define an invariant for unitary involutions with trivial discriminant algebra on split algebras. Nevertheless, this invariant does not extend in a functorial way to the non split case, see \S\ref{absolute.section}. Using the Rost invariant, one may still define a relative Arason invariant, satisfying formal properties established in~\S\ref{prop.section}. In particular, even though it is generally represented by cohomology classes of larger order, the relative invariant has order $2$, and satisfies some additivity property, see \cref{Chasles} and \cref{order}. 

An important feature compared  to previous papers is that we do not restrict to involutions with trivial lower-degree invariants, but also consider pairs of unitary involutions with isomorphic discriminant algebras, and pairs of orthogonal involutions with isomorphic Clifford algebras, see \S\ref{relative.section} and \ref{orthogonalinvariant.sec}. This leads to a notion of relative Arason invariant, with values in the quotient of $H^3(F)$ by the image under the Rost invariant of cocycles coming from the center of the relevant groups; this image is described for each type of group in~\cite{MPT} and~\cite{GM}. The relative Arason invariant has already proved to be useful in unitary type : Merkurjev considered it in the particular case of degree $4$ split algebras, and he used it to extend to groups of outer type $^2{\sf A}_3$ a theorem of Rost on $R$-equivalence classes, see~\cite[\S 4]{Me2000}. 
 In orthogonal type, an additional condition is required for defining the relative Arason invariant. This condition may be understood via the study of outer automorphisms of algebraic groups carried out in~\cite{QT-AutExt}, or in terms of similarity classes of hermitian forms, see Remarks~\ref{relorth.rem}(2) and~\ref{sim.rem}. It did not occur in previous papers, since only inner forms of the underlying groups were considered there. 

With those Arason invariants in hand, we can provide precise comparison theorems for degree 3 invariants of different type in several situations.  
Namely, when an algebra with $F'/F$-unitary involution $(B,\tau)$ is isomorphic to a tensor product $(C,\gamma)\otimes_F (F',\iota)$ for some $F$-algebra with symplectic or orthogonal involution $(C,\gamma)$, where $\iota$ denotes the unique non trivial automorphism of the quadratic extension $F'/F$, we call $(C,\gamma)$ a symplectic or an orthogonal descent of $(B,\tau)$. Conversely, if an algebra with symplectic or orthogonal involution over $F$ contains a subfield isomorphic to $F'$ with centralizer isomorphic to $(B,\tau)$ as an algebra with involution, we call it a quadratic extension of $(B,\tau)$. In both cases, strong relations exist between the degree $3$ invariants of the relevant types; see~\S\ref{symp.section} and~\ref{orth.section} for precise statements. 

Finally, we study properties of the Arason invariant for algebras of small degree with unitary involutions. In particular, we prove it is classifying in degree $4$. More precisely, two $F'/F$ unitary involutions are $F$-isomorphic if and only if their relative Arason invariant is trivial. In degree $4$ and $6$, we prove that hyperbolic involutions are characterized by the vanishing of the hyperbolic Arason invariant. 
This is not true anymore in degree $8$, where one would expect a characterization of totally decomposable involutions. For algebras of degree $8$ and index $4$, we use the relative Arason invariant to define a new absolute invariant, taking an arbitrary totally decomposable involution as a base point. We prove the value of this invariant does not depend on the chosen base point, using a cohomological invariant of the underlying biquaternion algebra introduced and studied in~\cite{Bar14}. This invariant does not coincide with the hyperbolic Arason invariant as, surprisingly, the hyperbolic involution is not totally decomposable in general, see \cref{hypversustotdec} and~\cref{hypnotdec}. 

This article owes a lot to Jean-Pierre Tignol, and all three of us are extremely grateful for what we learned from him. In particular, he supervised Alexandre Masquelein's thesis, whose main results are contained in this paper. 
The third named author also thanks Skip Garibaldi and Philippe Gille for useful conversations. Both authors thank the referee for his valuable comments.
 
\subsection*{Notations and preliminary results} 
\label{prel.section}
Let $F$ be a field. Throughout the paper, we assume F has characteristic different from 2. This may be unnecessary, as the Arason invariant, the Rost invariant, and the image of central cocycles are known in all characteristic. Nevertheless, we wish to avoid technical complications, notably when dealing with Arason invariant for groups of orthogonal type, for which the current literature makes this assumption. 

We use the notations and terminology of~\cite{KMRT} for algebras with involution. 
In particular, throughout the paper, we call $(B,\tau)$ an algebra with unitary involution over $F$ when $B$ is a finite dimensional $F$-algebra, with center a quadratic \'etale extension $F'=F[X]/(X^2-\delta)$ of $F$, such that $B$ is either simple (if $F'$ is a field) or a direct product of two simple algebras (if $F'\simeq F\times F$) and $\tau$ is an $F'/F$ semi-linear involution. In particular, $\tau$ acts on $F'$ as the unique non-trivial $F$-automorphism, which we denote by $\iota$. For any field extension $L/F$, the pair \[(B_L,\tau_L)=(B\otimes_F L,\tau\otimes \id)\] is an algebra with unitary involution over $L$, with center $L'=F'\otimes_F L$. 

The central simple algebra $B$ over $F'$ is endowed with an $F'/F$-unitary involution $\tau$ if and only if $B$ has split corestriction, or equivalently, its Brauer class has trivial corestriction, see~\cite[(3.1) \& (3.20)]{KMRT}. Under this condition, if $D$ is an $F'$-central division algebra Brauer equivalent to $B$, then $D$  is also endowed with an $F'/F$-unitary involution $\theta$. Moreover, there exists a hermitian module $(M,h)$ with values in $(D,\theta)$, such that $(B,\tau)\simeq(\End_D(M),\ad_h)$, where $\ad_h$ denotes the adjoint involution. The hermitian form $h$ is uniquely determined up to a scalar factor $\lambda\in F^\times$. 

If $B$ has even degree, $n=2m$, we let $\cd(\tau)$ be the discriminant algebra of $(B,\tau)$~\cite[\S10]{KMRT}. The discriminant algebra is a central simple algebra over $F$, and its Brauer class can be computed in some cases, as recalled below for further use. 
Assume $(B,\tau)$ admits an orthogonal descent, $(B,\tau)\simeq (A,\sigma)\otimes_F (F',\iota)$, for some central simple algebra with orthogonal involution $(A,\sigma)$ over $F$. Then by~\cite[(10.33)]{KMRT}, the Brauer class of the discriminant algebra is given by
\begin{equation}
\label{DiscOrthDescent.eq}
\cd(\tau)\sim\left\{\begin{array}{ll} 
\bigl(\delta,\disc(\sigma)\bigr)&\mbox{ if $m$ is even;}\\
A\otimes \bigl(\delta,\disc(\sigma)\bigr)&\mbox{ if $m$ is odd,}\\
\end{array}\right.
\end{equation}
where $\disc(\sigma)$ stands for the discriminant of the orthogonal involution $\sigma$. 
Moreover, a unitary involution on a quaternion algebra has a unique symplectic descent by~\cite[(2.22)]{KMRT}, that is 
$(Q,\tau)\simeq (H,\gamma)\otimes (F',\iota)$ for some uniquely determined quaternion algebra $H$ over $F$, with canonical (symplectic) involution $\gamma$. 
In view of~\cite[(10.30)]{KMRT}, we then have $\cd(\tau)\sim H$. Hence, two unitary involutions on a quaternion algebra are isomorphic if and only if they have isomorphic discriminant algebras. 
In particular, one may check using~\eqref{DiscOrthDescent.eq} that for all $a\in F^\times$, we have 
\begin{equation}
\label{QuatUnit.eq}
(M_2(F),\ad_\pform{a})\otimes_F (F',\iota)\simeq \bigl((\delta,a),\gamma\bigr)\otimes_F (F',\iota).
\end{equation}

Assume now $F'$ is a field and $B$ is split, that is $B\simeq\End_{F'}(V)$ for some $n$-dimensional vector space $V$ over $F'$. The involution $\tau\simeq\ad_h$ is adjoint to a hermitian form $h:\,V\times V\rightarrow (F',\iota)$, and again $h$ is unique up to a scalar factor $\lambda\in F^\times$. One may associate to $h$ a quadratic form called its Jacobson's trace and 
defined by $q_h(x)=h(x,x)\in F$ for all $x\in V$, see~\cite[Chap. 10 \S 1]{Scharlau}. 
So $q_h$ is a $2n$-dimensional quadratic form over $F$, and by {\em loc. cit.} Thm 1.1, the isometry class of $q_h$ determines the isometry class of $h$. 
If $h$ diagonalises as $h\simeq\qform{a_1,\dots,a_n}$ for some $a_i\in F^\times$, then 
$q_h\simeq\qform{1,-\delta}\otimes  \qform{a_1,\dots,a_n}$. 
The invariants of $q_h$ are computed in~\cite[p. 350]{Scharlau}. Its discriminant and full Clifford algebra satisfy 
\begin{multline}
\label{discJacobson}
d(q_h)=\left\{\begin{array}{ll} 
\delta\in F^\times/F^{\times 2}&\mbox{ if $n=2m+1$ is odd,}\\
1\in F^\times/F^{\times 2}&\mbox{ if $n=2m$ is even.}\\
\end{array}\right.\\
\mbox{ and }\C(q_h)\sim\bigr(\delta,d(h)\bigl)
\mbox{ where }d(h)=(-1)^\frac{n(n-1)}{2}a_1\dots a_n\in F^\times/N_{F'/F}(F'^\times)\\ 
\mbox{ is the discriminant of the hermitian form $h$.}\\
\end{multline}
In particular, if $n=2m$ is even, by the structure theorem for Clifford algebras, see~\cite[V. Thm. 2.5(3)]{Lam}, we have $\C_0(q_h)\simeq\C_+\times\C_-$ with $\C_+\sim\C_-\sim\C(q_h)$. Comparing with~\cite[(10.35)]{KMRT}, we get that $\C_+$ and $\C_-$ are Brauer equivalent to $\cd(\ad_h)$. 
Assume now that $n=2m+1$ is odd. The even Clifford algebra $\C_0(q_h)$ is the centraliser of $F'$ in $\C(q_h)$, hence, $\C_0(q_h)\sim \C(q_h)_{F'}$, see~\cite[Chap. 8, Thm. 5.4]{Scharlau}. In view of the formula above for $\C(q_h)$, we have $\C(q_h)_{F'}\sim 0$. So we finally get 
\begin{multline}
\label{CliffordJacobson}
\C_0(q_h)\simeq\left\{\begin{array}{ll} 
0\in\br(F')&\mbox{ if $n=2m+1$ is odd,}\\
\C_+\times\C_-\in\br(F)\times\br(F),&\mbox{ if $n=2m$ is even,}\\ 
\end{array}\right.\\
\mbox{ with }\C_+\sim \C_-\sim\cd(\ad_h).\\
\end{multline}

For any discrete torsion Galois module $M$, we let $H^i(F,M)$ be the Galois cohomology group $H^i\bigl({\mathrm {Gal}}(F_{\mathrm {sep}}/F),M\bigr)$.  If $G$ is a group scheme over $F$, ${\mathrm {Gal}}(F_{\mathrm {sep}}/F)$ acts continuously on $G(F_{\mathrm {sep}})$ and we let $H^1(F,G)=H^1\bigl({\mathrm {Gal}}(F_{\mathrm {sep}}/F),G(F_{\mathrm {sep}})\bigr)$. We also use the notation \[H^i(F)=H^i\bigl(F,\mq/\mz(i-1)\bigr)\mbox{, see~\cite[Appendix A, p. 151]{GMS}.}\]
In all three cases, we denote by $\res_{F/L}$, respectively $\cores_{L/F}$, the restriction and corestriction maps in Galois cohomology, where $L/F$ is an arbitrary (respectively a finite degree) field extension. For each integer $m\geq 0$, $_mH^i(F)$ denotes the $m$-torsion subgroup of $H^i(F)$. Using the norm residue isomorphism, one may check that $_2H^i(F)=H^i(F,\mu_2)$ and $_4H^3(F)=H^3(F,\mu_4^{\otimes 2})$. In particular, $_2H^1(F)\simeq F^\times/F^{\times 2}$ and we let $(a)\in \,_2H^1(F)$ be the square class of $a\in F^\times$ and 
$(a_1,\dots,a_i)=(a_1)\cdot\dots\cdot (a_i)\in \,_2H^i(F)$ the cup-product of the square classes of $a_1,\dots, a_i\in F^\times$. In particular, we use the same notation $(a_1,a_2)$ for the quaternion $F$-algebra and its Brauer class in $_2\br(F)\simeq\, _2H^2(F)$. 
Recall that the Milnor $K$-ring $K_* (F)$ acts on $\sqcup _{i\geq 0}H^i(F)$.
In particular, for any $\alpha\in\br(F)\simeq H^2(F)$,  $F^\times\cdot\alpha\subset H^3(F)$ is a subgroup, and we let 
\[M^3_\alpha(F)=H^3(F)/(F^\times\cdot\alpha).\]
Given in addition a quadratic extension $F'/F$ and another Brauer class $\beta\in\br F'$, we define in a similar way
\[N^3_{\alpha,\beta}(F)=H^3(F)/\bigr(F^\times\cdot\alpha+\cores_{F'/F}(F'^\times\cdot\beta)\bigl).\]
In particular, we have $N^3_{\alpha,0}(F)=M^3_\alpha(F)$. 

We denote by $I(F)$ the fundamental ideal in the Witt ring $W(F)$. For any (non degenerate) quadratic form $q$ over $F$, if the Witt class of $q$ is in $I^3(F)$, we let $e_3(q)\in\, _2H^3(F)$ be its Arason invariant. 
If $G/F$ is an absolutely almost simple simply connected algebraic group, we let $\rho_G$ be the corresponding Rost invariant, \[\rho_G:\, H^1(F,G)\rightarrow H^3(F)\mbox{, see~\cite{GMS}.}\] 
An explicit description of $\rho_G$ is given in some cases in~\cite[p. 107--108]{GMS}. In particular, if $G=\Spin(q_0)$ is the Spin group of some non degenerate quadratic form $q_0$, then $\rho_G$ can be described in terms of the Arason invariant.  

Let $i:\,G\rightarrow G'$ be a homomorphism of absolutely almost simple simply connected algebraic groups. By~\cite[Thm. 9.11]{GMS}, 
 for all field extension $L/F$, we have 
 \[\rho_{G'_L}(i^{(1)}(\xi))=n_i \,\rho_{G_L}(\xi)\in H^3(L)\mbox{ for all }\xi\in H^1(L,G_L),
\]
where $i^{(1)}:\, H^1(L,G_L)\rightarrow H^1(L, G'_L)$ denotes the map induced by $i$, and the integer $n_i$ is the index defined by Dynkin in~\cite{Dynkin}, also called the Rost multiplier of $i$.
We will use the following : 
\begin{lem}
\label{rostmultiplier}
The homomorphism 
\[j:\,\SL_n(F)\rightarrow \SL_{2n}(F), \ \ X\mapsto \begin{pmatrix} X&0\\0&(X^{-1})^{t}\end{pmatrix}\]
has Rost multiplier $n_j=2$. 
\end{lem}
\begin{proof} 
The homomorphism $j$ is the direct sum of the tautological representation $\rho$ of $\SL_n$ and its dual representation $\rho^*$. By definition of the index, we have $n_\rho=1$. Moreover, composition with a group automorphism does not change the index, so we also have $n_{\rho^*}=1$. 
The result follows by~\cite[Prop. 7.9(5)]{GMS}. 
\end{proof} 

Our main interest in this paper is for outer groups of type $^2\mathrm{A_{n-1}}$. So, let us assume $G=\SU(B,\tau_0)$ for  some degree $n$ algebra with unitary $F'/F$-involution $(B,\tau_0)$. 
We recall from~\cite[(29.18)]{KMRT} that the cohomology set $H^1\bigl(F,\SU(B,\tau_0)\bigr)$ may be identified with $\SSym(B,\tau_0)^\times/\approx$, where 
\begin{equation}
\label{ssym.eq}
\SSym(B,\tau_0)=\{(s,z)\in\Sym(B,\tau_0)^\times\times F'^\times,\ \Nrd_B(s)=N_{F'/F}(z)\},
\end{equation}
and the equivalence relation $\approx$ is defined by $(s,z)\approx(s',z')$ if and only if there exists $b\in B^\times$ such that $s'=bs\tau(b)$ and $z'=\Nrd_B(b) z$. 
Let $R_{F'/F}(\mu_n)$ be the corestriction of the algebraic group $\mu_n$ of roots of unity. 
The center of $G$ is isomorphic to the kernel $\mu_{n[F']}$ of the norm map 
$N_{F'/F}:\,R_{F'/F}(\mu_n)\rightarrow \mu_n.$
Under~\eqref{ssym.eq}, $H^1(F,\mu_{n[F']})$ corresponds  to the classes of all pairs $(x,z)\in F^\times\times F'^\times$ such that $x^n=N_{F'/F}(z)$. 
From the description of the image under $\rho_G$ of the class of such an element given in~\cite[Thm. 1.10 \& Thm. 1.11]{MPT} and ~\cite[\S 10]{GM}, we get 
\begin{equation}
\label{RostCenter.eq}
\rho_{\SU(B,\tau_0)}\bigl(H^1(F,\mu_{n[F']})\bigr)\subset\left\{\begin{array}{ll}
\cores_{F'/F}(F'^\times\cdot[B])&\mbox{ if $n$ is odd,}\\
F^\times\cdot[\cd(\tau_0)]+\cores_{F'/F}(F'^\times\cdot[B])&\mbox{ if $n$ is even.}\\
\end{array}\right.
\end{equation}

\section{A degree $3$ invariant for unitary involutions} 
\label{def.section}
The Arason invariant of quadratic forms may be used to define a degree $3$ invariant for unitary involutions on split algebras, provided the underlying algebra $B=\End_{F'}(V)$ has even degree $n$ and the involution has a split discriminant algebra, as we shall explain in \S\ref{absolute.section}. This invariant does not extend in a functorial way to the non-split case. Nevertheless, using the Rost invariant, we may still define a relative Arason invariant for unitary involutions, see~\S\ref{relative.section}. Choosing a particular base point, this produces absolute Arason invariants, under some additional condition on the underlying algebra. Notably, if the algebra has even co-index, so that it admits a hyperbolic involution, which is unique up to isomorphism, we get a so-called {\em hyperbolic} Arason invariant, see~\S\ref{hyperbolic.section}. Finally, multiplying the Rost invariant by $2$, we define in~\S\ref{f3.section} a relative invariant with values in $H^3(F)$ when $B$ has exponent at most $2$. 

From now on, we let $F'$ be a quadratic \'etale extension of $F$, $F'=F[X]/(X^2-\delta)$ for some $\delta\in F^\times$, and we denote by $\iota$ the unique non-trivial $F$-automorphism of $F'$.  

\subsection{Absolute degree $3$ invariant in the split case}
 \label{absolute.section}
 
Let $(B,\tau)$ be an algebra of degree $n$ with $F'/F$-unitary involution such that the algebra $B$ is split, and consider a hermitian space $(V,h)$ over $(F',\iota)$ such that $B\simeq\End_{F'}(V)$ and $\tau\simeq\ad_h$. We assume in addition that $n=2m$ is even and the discriminant algebra $\cd(\tau)$ is split. As recalled in the preliminary section, it follows the Jacobson's trace $q_h$ of $h$ belongs to $I^3(F)$. Hence, one may consider the Arason invariant $e_3(q_h)\in\, _2H^3(F)$. Moreover, for all $\lambda\in F^\times$, we have $q_{\qform{\lambda}h}=\qform{\lambda}q_h$ and $e_3(\qform{\lambda}q_h)=e_3(q_h)\in\, _2H^3(F)$. Hence, the cohomology class $e_3(q_h)$ is a well defined invariant of the algebra with involution $(\End_{F'}(V),\ad_h)$. We will refer to it as its Arason invariant, and use the notation \[e_3(\ad_h)=e_3(q_h)\in\,_2H^3(F).\] This invariant clearly is functorial, that is for any field extension $L/F$, we have $e_3\bigl((\ad_h)_L\bigr)=\res_{L/F}\bigr(e_3(\ad_h)\bigl)\in\, _2H^3(L)$.

As for orthogonal involutions (see~\cite[\S 3.4]{BPQ}), we claim this invariant does not extend in a functorial way to the non-split case, at least when the degree $n$ of $A$ satisfies $n\equiv 4\mod 8$. Specifically, fix a non-zero integer $n$ such that $n\equiv 4\mod 8$, and assume for the sake of contradiction that there exists a degree $3$ invariant satisfying the following conditions: 
\begin{enumerate}
\item[(i)] $e_3(\tau)$ is defined for all degree $n$ algebra with unitary involution $(B,\tau)$ over a field extension $L$ of the base field $F$ such that the discriminant algebra $\cd(\tau)$ is split; 
\item[(ii)] $e_3(\tau)$ has values in $H^3(L)/H_B$, where $H_B$ is a subgroup of $H^3(L)$, depending on the algebra $B$;
\item[(iii)] If $B$ is split and $\tau$ is adjoint to a hermitian form $h$, then the group $H_B$ is trivial and $e_3(\tau)=e_3(\ad_h)=e_3(q_h)\in\,_2 H^3(L)$; 
\item[(iv)] Given a field extension $E$ of $L$, the restriction map in cohomology induces a map $\res_{E/L}:\, H^3(L)/H_B\rightarrow H^3(E)/H_{B_{E}}$ and we have \[e_3(\tau_{E})=\res_{E/L}\bigl(e_3(\tau)\bigr)\in H^3(E)/H_{B_{E}}.\]
\end{enumerate}
Assume first that $n=4$. 
In order to get a contradiction, consider a decomposable degree $4$ algebra with unitary involution $(B,\tau)$, that is 
\[(B,\tau)=(Q_1,\gamma_1)\otimes (Q_2,\gamma_2)\otimes (L',\iota), \]
for some quaternion algebras $Q_i$ with center $L$ and canonical involution $\gamma_i$ for $i=1,2$, and a quadratic field extension $L'/L$. 
For a suitable choice of $L$, $Q_1$, $Q_2$ and $L'$, we may assume $B$ is division. 
Moreover, $\gamma_1\otimes\gamma_2$ is an orthogonal involution with trivial discriminant by~\cite[(7.3)(4)]{KMRT},  hence the discriminant algebra $\cd(\tau)$ is split, see~\eqref{DiscOrthDescent.eq}. 
Therefore, we may consider $e_3(\tau)\in H^3(L)/H_B$. Pick an element $\xi\in H^3(L)$ such that $e_3(\tau)=\xi\mod H_B$. We claim $\xi\in \,_4H^3(L)$. Indeed, there exists a biquadratic field extension $E/L$, which is a subfield of $Q_1\otimes_FQ_2$, such that both $Q_1$ and $Q_2$ are split over $E$. Hence, the algebra with involution $(B_{E},\tau_{E})$ is split and hyperbolic. It follows $\tau_{E}$ is adjoint to a hyperbolic hermitian form $h$, and its Jacobson's trace $q_h$ is also hyperbolic. Therefore $H_{B_{E}}$ is trivial and  $e_3(\tau_{E})=0\in H^3(E)$ by condition (iii). Using condition (iv), we get that $\res_{E/L}(\xi)=0\in H^3(E)$ and applying the corestriction map in cohomology, this leads to $4\xi=0$ as claimed.

Now, since $\xi$ is a sum of $2$-power order cohomology classes in $H^3(F)$, there exists a family of $3$-fold Pfister forms $\pi_1,\dots,\pi_r$ such that $\xi$ vanishes over the function field $\cf$ of the product of the corresponding quadrics, that is $\res_{\cf/L}(\xi)=0$. In particular, applying again (iv), we have $e_3(\tau_\cf)=0$. On the other hand, by Merkurjev's index reduction formula~\cite{Me:IR} (see also~\cite{Kahn}), $B_\cf$ is still division, so that the involution $\tau_\cf$ is anisotropic. This leads to a contradiction as follows. Let $\cf_1$ be the function field of the Weil restriction $R_{\cf'/\cf}(\SB_{B_\cf})$ of the Severi-Brauer variety of $B_\cf$, where $\cf'=L'\otimes_L\cf$ is the center of $B_\cf$. By (iv), we have $e_3(\tau_{\cf_1})=\res_{\cf_1/\cf}(e_3(\tau_\cf))=0$. In addition, the field $\cf_1$ is a generic splitting field of $B_\cf$, therefore $B_{\cf_1}$ is split and $\tau_{\cf_1}$ is the adjoint of a $4$-dimensional hermitian form $h_\tau$, whose trace form $q_{h_\tau}$ is an $8$-dimensional quadratic form. It follows by (iii) that $q_{h_\tau}$ is hyperbolic, that is $h_\tau$ and $\tau_{\cf_1}$ are also hyperbolic. This contradicts~\cite[Thm. A.2]{K:isotropy}, which asserts that the anisotropic involution $\tau_\cf$ remains anisotropic over $\cf_1$. Hence, there is no invariant satisfying conditions (i) to (iv) above with $n=4$. 

Assume now $n=4+8k$ for some integer $k$, and consider the algebra $M_{2k+1}(B)$, with $B$ as above, endowed with the involution adjoint to the hermitian form with values in $(B,\tau)$ defined by $h=\qform{1}+k{\mathbb H}$, where ${\mathbb H}$ denotes a hyperbolic plane. Since the anisotropic kernel of this form is $(B,\tau)$, it has the same invariants as $(B,\tau)$. Therefore, the same argument as above leads to a contradiction, and this concludes the proof.

\begin{remark}
\label{deg3.1} 
If $B$ has odd degree, the quadratic form $q_h$ has discriminant $\delta$, hence it does not belong to $I^3(F)$. 
When $\deg(B)=3$, we may nevertheless define an absolute invariant for unitary involutions, even when $B$ is non-split. This is explained in~\cite[\S 19.B]{KMRT} : to any unitary involution $\tau$ of $B$, corresponds a uniquely defined $3$-fold Pfister form $\pi(\tau)$, so that $e_3\bigl(\pi(\tau)\bigr)\in H^3(F,\mu_2)$ is a well defined invariant of $\tau$, see~\cite[(30.21)]{KMRT}. It follows from the definition of $\pi(\tau)$ that this invariant is functorial. Moreover, it  is classifying by~\cite[(19.6)]{KMRT} : two unitary involutions $\tau_0$ and $\tau$ of $B$ are isomorphic if and only if 
$e_3(\pi(\tau_0))=e_3(\pi(\tau))\in H^3(F,\mu_2)$. 
\end{remark}

\subsection{Relative degree $3$ invariant}
\label{relative.section}

Using the Rost invariant, one may still define a relative invariant for unitary involutions. This definition is already contained in~\cite[\S 4.3]{Tignol:hyderabad} and~\cite{M:thesis}, under the more restrictive assumption that both involutions have split discriminant algebra when $B$ has even degree. 

Consider two $F'/F$-unitary involutions $\tau_0$ and $\tau$ of the algebra $B$, and assume in addition that $\tau_0$ and $\tau$ have isomorphic discriminant algebras when $B$ has even degree. As explained in~\cite[\S 29.D]{KMRT}, we may associate to $\tau$ a cohomology class $\eta\in H^1\bigl(F,\PGU(B,\tau_0)\bigr)$ representing the triple $[B,\tau,\id_{F'}]$. The exact sequence 
\begin{equation}
\label{ExactSeq.eq}
1\rightarrow\mu_{n[F']}\rightarrow \SU(B,\tau_0)\rightarrow\PGU(B,\tau_0)\rightarrow 1
\end{equation}
induces a connecting map 
\[\partial:\,H^1(F,\PGU(B,\tau_0))\rightarrow H^2(F,\mu_{n[F']}). \]
By a twisting argument as in~\cite[I.5.4]{Serre} (see also~\cite[\S 1]{Gar}), one may check that $\partial(\eta)$ is the difference between the Tits classes of the groups $\SU(B,\tau)$ and $\SU(B,\tau_0)$. 
If $B$ has odd degree, $n=2m+1$, the restriction map identifies $H^2(F,\mu_{n[F']})$ with the kernel of the corestriction $\cores_{F'/F}:\,_n\br(F')\rightarrow\, _n\br(F)$, see~\cite[Rem.  p. 309]{CTGP}. By~\cite[(31.8)]{KMRT}, the groups $\SU(B,\tau_0)$ and $\SU(B,\tau)$ both have the same Tits class, which corresponds to $[B]$ under this identification. 
Assume now $B$ has even degree $n=2m$. We may identify $H^2(F,\mu_{n[F']})$ with a subgroup of $\br(F)\times\br(F')$ as in~\cite[Prop. 2.10]{CTGP}. The Tits classes of $\SU(B,\tau_0)$ and $\SU(B,\tau)$ respectively correspond to the pairs $([\cd(\tau_0)],[B])$ and $([\cd(\tau)],[B])$ by~\cite[(31.8)]{KMRT}. Hence in both cases, we get $\partial(\eta)=0$ so that $\eta$ lifts to a class $\xi\in H^1(F,\SU(B,\tau_0))$, which is uniquely defined up to the action of $H^1(F,\mu_{n[F']})$. 
In view of the behavior of the Rost invariant under twisting, as described in~\cite[Lemme 7]{Gille} (see also~\cite[Prop. 1.7]{MPT}), it follows that $\rho_{\SU(B,\tau_0)}(\xi)\in H^3(F)$ is well defined up to an element which belongs to the image of $H^1(F,\mu_{n[F']})$ under $\rho_{\SU(B,\tau_0)}$, hence to the subgroup of $H^3(F)$ given in~\eqref{RostCenter.eq}.

This leads to the following : 
\begin{defi} 
Given two $F'/F$-unitary involutions $\tau_0$ and $\tau$ of the degree $n$ algebra $B$, 
with $\cd(\tau_0)\sim\cd(\tau)$ if $n$ is even, we define their relative Arason invariant by 
\[e_3^{\tau_0}(\tau)=\rho_{\SU(B,\tau_0)}(\xi)\in N^3_{\alpha,\beta}(F),\]
where 
$\xi\in H^1\bigr(F,\SU(B,\tau_0)\bigl)$ is a cocycle with image in 
$H^1\bigl(F,\PGU(B,\tau_0)\bigr)$ corresponding to the isomorphism class of $[B,\tau,\id_{F'}]$, 
\[\alpha=\left\{\begin{array}{ll}
[\cd(\tau_0)]\in\br(F)&\mbox{if $n$ is even,}\\
0&\mbox{if $n$ is odd,}\\
\end{array}\right.\mbox{ and }\beta=[B]\in\br(F').\]
\end{defi}

We will refer to $e_3^{\tau_0}$ as the relative Arason invariant, with $\tau_0$ as a base point. Depending on the degree $n$ of $B$, it has values in   
\[\left\{\begin{array}{l}
H^3(F,\mq/\mz(2))/\cores_{F'/F}(F'^\times\cdot [B])\mbox{ if $n$ is odd, and}\\
H^3(F,\mq/\mz(2))/\bigl(F^\times\cdot[\cd(\tau_0)]+\cores_{F'/F}(F'^\times\cdot [B])\bigr)\mbox{ if $n$ is even.}
\end{array}\right.\] 
If $\tau$ is isomorphic to $\tau_0$, then $[B,\tau,\id_{F'}]$ corresponds to the base point in the cohomology set $H^1\bigl(F,\PGU(B,\tau_0)\bigr)$, hence $e_3^{\tau_0}(\tau)=0$ in this case. 
The invariant $e_3^{\tau_0}$ is a functorial invariant. Indeed, for all field extension $L/F$, the restriction map in Galois cohomology induces a map 
\begin{equation}
\res_{L/F}:\,N^3_{\alpha,\beta}(F)\rightarrow N^3_{\alpha_L,\beta_L}(L),
\end{equation}
where $\alpha_L=\res_{L/F}(\alpha)$, and 
$\beta_L=[B\otimes_F L]\in \br(L')\mbox{ with }L'=F'\otimes_F L.$
In particular, if $n$ is even, we have $\alpha_L=[\cd(\tau_0)_L]=[\cd({\tau_0}_L)]\in\br(L)$. 
Moreover, in both cases, \[e_3^{{\tau_0}_L}(\tau_L)=\res_{L/F}\bigr(e_3^{\tau_0}(\tau)\bigl)\in N^3_{\alpha_L,\beta_L}(L) .\]

\begin{example}
\label{deg3.2}
If $B$ has degree $3$, then the relative Arason invariant classifies unitary involutions up to isomorphism. 
This follows from the description of the Rost invariant given in~\cite[(31.45)]{KMRT}. 
Indeed, as explained there, $\rho_{\SU(B,\tau_0)}(\xi)$ belongs to $H^3(F,\mu_6^{\otimes 2})$, and its $2$-primary part is equal to $e_3\bigl(\pi(\tau)\bigr)-e_3\bigl(\pi(\tau_0)\bigr)$. On the other hand, $\cores_{F'/F}\bigl(F'^\times\cdot [B]\bigr)$ is contained in $H^3(F,\mu_3^{\otimes 2})$. 
Therefore, if $e_3^{\tau_0}(\tau)=0\in H^3(F)/\cores_{F'/F}\bigl(F'^\times\cdot [B]\bigr)$, the $2$-primary part of $\rho_{\SU(B,\tau_0)}(\xi)$ is trivial, that is $e_3\bigl(\pi(\tau)\bigr)=e_3\bigl(\pi(\tau_0)\bigr)\in H^3(F,\mu_2)$. This implies $\tau\simeq \tau_0$ by~\cite[(19.6)]{KMRT}, see also Remark~\ref{deg3.1}. 
\end{example}

In the split case, the invariant $e_3^{\tau_0}$ can be described in terms of the the Arason invariant of quadratic forms as follows. 
\begin{lem}
\label{split}
Let $V$ be a vector space of dimension $n$ over $F'$, and $h_0$ and $h$ two hermitian forms on $V$, with the same discriminant $d(h_0)=d(h)\in F^\times/N_{F'/F}(F'^\times)$. 
The involutions $\tau_0=\ad_{h_0}$ and $\tau=\ad_h$ have isomorphic discriminant algebras when $n$ is even, and their relative Arason invariant is given by 
\[e_3^{\tau_0}(\tau)=e_3(q_h-q_{h_0})\in M^3_{\alpha}(F), \]
where $q_h$ and $q_{h_0}$ are Jacobson's traces of the hermitian forms $h$ and $h_0$, and $\alpha$ is $0$ if $n$ is odd and $[\cd(\tau_0)]\in\br(F)$ if $n$ is even. 
\end{lem} 

\begin{proof}
First of all, since $B$ is split, $\beta=[B]=0\in\br(F')$ and $e_3^{\tau_0}$ has values in $N^3_{\alpha,0}(F)=M^3_\alpha(F)$. 
Pick a diagonalisation for the hermitian forms, $h_0=\qform{a_1,\dots, a_n}$ and $h=\qform{b_1,\dots, b_n}$ where $a_i,b_i\in F^\times$ for $1\leq i\leq n$. By assumption, there exists $y\in F'^\times$ such that $a_1\dots a_n =b_1\dots b_n N_{F'/F}(y)$. Therefore, we have \[q_h-q_{h_0}=\qform{1,-\delta}\qform{b_1,\dots, b_n,-a_1,\dots, -a_n}=q\perp\qform{-\delta} q,\] with $q$ of dimension $2n$ and discriminant $N_{F'/F}(y)$. So $q_h-q_{h_0}$ has trivial discriminant. Moreover, its Clifford invariant can be computed using~\cite[V (3.15)\&(3.16)]{Lam}
\[e_2(q_h-q_{h_0})=e_2(q)+e_2(-\delta N_{F'/F}(y) q)=e_2(q)+e_2(q)+(-\delta N_{F'/F}(y),N_{F'/F}(y))=0.\]
Hence, under the conditions of the lemma, the difference $q_h-q_{h_0}$ belongs to $I^3(F)$, and we may consider $e_3(q_h-q_{h_0})\in H^3(F)$. It remains to prove its image in $M^3_{\alpha}(F)$ coincides with $e_3^{\tau_0}(\tau)$. 

By~\cite[(29.19)]{KMRT}, the isometry class of $h$ corresponds to a cohomology class 
\[\xi\in H^1(F,\SU(V,h_0))=H^1(F,\SU(B,\tau_0)).\] 
Under~\eqref{ssym.eq}, $\xi$ corresponds to the class of a pair $(s,z)$ where $s\in \End_{F'}(V)$ is a $\tau_0$-symmetric element and $z\in F'^\times$ satisfies $N_{F'/F}(z)={\mathrm{det}}(s)$. Moreover, we have $h\simeq h_s$, where $h_s$ is the hermitian form defined by $h_s(x,y)=h_0(s^{-1}(x),y)$. 
Since $\tau\simeq\ad_{h_s}\simeq \Int(s)\circ \tau_0$, the image of $\xi$ in $H^1(F,\GU(B,\tau_0))$ corresponds to the conjugacy class of $\tau$ under the bijection described in~\cite[(29.16)]{KMRT}. 
Hence, $\xi$ maps to a cohomology class corresponding to $[B,\tau,\id_{F'}]$ in $H^1(F,\PGU(B,\tau_0))$. 
By definition of the relative Arason invariant, we get
\[e_3^{\tau_0}(\tau)=\rho_{\SU(B,\tau_0)}(\xi)\in M^3_{\alpha}(F). \]
On the other hand, the  description of $\rho_{\SU(B,\tau_0)}$ given in~\cite[(31.44)]{KMRT} gives 
\[\rho_{\SU(B,\tau_0)}(\xi)=e_3(q_h-q_{h_0})\in H^3(F),\] and this concludes the proof. 
\end{proof}

\begin{remark} 
\label{splitcase.remark}
Given two involutions $\tau_0$ and $\tau$ on the split algebra $\End_{F'}(V)$, there exists hermitian forms $h_0$ and $h$ on $V$ such that $\tau\simeq\ad_h$ and $\tau_0\simeq\ad_{h_0}$. Those hermitian forms are uniquely defined up to a scalar factor. 

(i) If $n$ is odd, we may choose $h$ in its similarity class so that $h$ and $h_0$ have the same discriminant. The above lemma applied to this particular choice for $h$ computes $e_3^{\tau_0}(\tau)\in H^3(F)$. 

(ii) Assume now $n$ is even and $\tau_0$ and $\tau$ have isomorphic discriminant algebras. By~\cite[(10.35)]{KMRT}, the underlying hermitian forms $h_0$ and $h$ have the same discriminant in $F^\times/N_{F'/F}(F'^\times)$. We claim that the value of $e_3(q_h-q_{h_0})\in M_\alpha^3(F)$ does not depend on the choice of $h_0$ and $h$ in their respective similarity classes. 
Indeed, given $\lambda_0,\lambda\in F^\times$, we have 
\begin{align*}e_3(q_{\qform{\lambda}h}-q_{\qform{\lambda_0}h_0})=e_3(\qform{\lambda}q_h-\qform{\lambda_0}q_{h_0})\\
=e_3\bigl(\qform{\lambda}(q_h-q_{h_0})+\qform{\lambda}\pform{\lambda^{-1}\lambda_0}q_{h_0}\bigr).
\end{align*}
Since both $q_h-q_{h_0}$ and $\pform{\lambda^{-1}\lambda_0}q_{h_0}$ are in $I^3(F)$, we get 
\begin{align*}
e_3(q_{\qform{\lambda}h}-q_{\qform{\lambda_0}h_0})=e_3\bigl(\qform{\lambda}(q_h-q_{h_0})\bigr)+e_3\bigl(\qform{\lambda}\pform{\lambda^{-1}\lambda_0}q_{h_0}\bigr)\\
=e_3(q_h-q_{h_0})+(\lambda^{-1}\lambda_0)\cdot [\cd(\tau_0)],
\end{align*}
where the last equality follows from the fact that the Clifford invariant of $q_{h_0}$ is $[\cd(\tau_0)]$. 
Hence, replacing $h$ and $h_0$ by a scalar multiple does not change the class of $e_3(q_h-q_{h_0})\in H^3(F)/F^\times\cdot [\cd(\tau_0)],$ as required.  

(iii) Finally, assume in addition that the discriminant algebras $\cd(\tau)$ and $\cd(\tau_0)$ are split, so that $e_3^{\tau_0}$ has values in $H^3(F)$. Then both $q_h$ and $q_{h_0}$ belong to $I^3(F)$. Hence we have 
\[e_3^{\tau_0}(\tau)=e_3(q_h-q_{h_0})=e_3(q_h)-e_3(q_{h_0})=e_3(\tau)-e_3(\tau_0)\in H^3(F).\]
In other words, the relative invariant coincides in this case with the difference of the absolute invariants of both involutions. 
\end{remark}

\subsection{Hyperbolic degree $3$ invariant}
\label{hyperbolic.section}
In this section, we assume that the algebra $B$ has even co-index, that is $B\simeq M_{2r}(D)$ for some integer $r\geq 1$ and some central division algebra $D$ over $F'$. In particular, this implies the degree $n$ of $B$ is even. We also assume $B$ admits $F'/F$-unitary involutions, hence it may be endowed with a hyperbolic $F'/F$-unitary involution $\tau_0$, which is unique up to isomorphism. 
As explained in~\cite[Proof of (10.36)]{KMRT}, the Brauer class of the discriminant algebra of $\tau_0$ can be computed after extending scalars to a generic splitting field of $B$; over such a field, $\tau_0$ is adjoint to a hyperbolic hermitian form, which has trivial discriminant, hence $\cd(\tau_0)$ is split. 

Let $\tau$ be a unitary $F'/F$ involution on $B\simeq M_{2r}(D)$ with split discriminant algebra. The hyperbolic Arason invariant of $\tau$ is defined by 
\[e_3^{\mathrm{hyp}}(\tau)=e_3^{\tau_0}(\tau)\in H^3(F)/\cores_{F'/F}(F'^\times\cdot[B]),\mbox{ where $\tau_0$ is hyperbolic}.\]
The invariant $e_3^{\mathrm{hyp}}$ is functorial, and it vanishes for hyperbolic involutions. 

\begin{example}
\label{splithyp.example}
(See~\cite[\S 4.1]{Tignol:hyderabad}). Assume $B$ is split and has even degree $n = 2m$. The hyperbolic degree $3$ invariant coincides in this case with the absolute invariant introduced in \S\ref{absolute.section}. This follows from~\cref{splitcase.remark}(iii). Indeed, consider two hermitian forms $h_0$ and $h$ over $(F',\iota)$ such that $\tau_0=\ad_{h_0}$ and $\tau=\ad_h$. Since $\tau_0$ is hyperbolic, $h_0$ and $q_{h_0}$ also are. Therefore, $e_3(\tau_0)=e_3(q_{h_0})=0$ and we get $e_3^{\tau_0}(\tau)=e_3(\tau)=e_3(q_h)$. 
\end{example}

\subsection{The $f_3^{\tau_0}$ invariant.}
\label{f3.section}
In this section, we assume that the algebra $B$ has exponent dividing $2$ and is endowed with a unitary involution $\tau_0$. It follows that the subgroup $\cores_{F'/F}(F'^\times\cdot [B])$ of $H^3(F)$ consists of classes that are of order at most $2$. 
Besides, if $n$ is even, the discriminant algebra $\cd(\tau_0)$ also has exponent $2$. Hence, by~\eqref{RostCenter.eq}, we have 
\begin{equation}
\label{2RostCenter.eq}
2\rho_{\SU(B,\tau_0)}\bigl(H^1(F,\mu_{n[F']})\bigr)=0.
\end{equation}
Therefore, the same procedure as in \S\ref{relative.section}, replacing the Rost invariant $\rho_{\SU(B,\tau_0)}$ by 
$2\rho_{\SU(B,\tau_0)}$, provides an invariant denoted by $f_3^{\tau_0}$ and with values in $H^3(F)$. 
More precisely, given two $F'/F$-unitary involutions $\tau_0$ and $\tau$ of $B$, and assuming $\cd(\tau_0)\sim\cd(\tau)$ if $n$ is even, we define their relative $f_3$ invariant by 
\begin{equation}
\label{f3.eq}
f_3^{\tau_0}(\tau)=2\rho_{\SU(B,\tau_0)}(\xi)\in H^3(F),
\end{equation}
where $\xi\in H^1(F,\SU(B,\tau_0))$ is any element with image in $H^1(F,\PGU(B,\tau_0))$ corresponding to the isomorphism class of $[B,\tau,\id_{F'}]$. 
We thus have $f_3^{\tau_0}(\tau)=2c$, where $c\in H^3(F)$ is an arbitrary cohomology class such that $e_3^{\tau_0}(\tau)$ is the class of $c$ in $N_{\alpha,\beta}^3(F)$. This is a well defined invariant since any two possible values of $c$ differ by an element of $F^\times\cdot\alpha+\cores_{F'/F}(F'^\times\cdot\beta)$, which is of order at most $2$. 

\begin{example}
If $B$ is split, then $f_3^{\tau_0}$ is identically zero. This follows from~\cref{split}, since the Arason invariant for quadratic forms has values in $_2H^3(F)$. It is also a consequence of the fact that the Rost invariant $\rho_{\SU(B,\tau_0)}$ has order $2$ in this case, see~\cite[Thm. 12.6]{GMS}. 
\end{example} 

\begin{example}
If $B$ has even co-index and exponent dividing $2$, for all unitary involution $\tau$ with split discriminant algebra, we define $f_3^{\hyp}(\tau)=f_3^{\tau_0}(\tau)$, where $\tau_0$ is a hyperbolic unitary involution of $B$. 
\end{example} 

\begin{example}
If $B$ has exponent $2$, the Rost invariant has order $4$ in general, see~\cite[Thm. 12.6]{GMS}. Therefore, there are examples of $(B,\tau_0,\tau)$ where $f_3^{\tau_0}(\tau)$ is non trivial. 
% Produire un exemple explicite. En degré 4 et algèbre discriminante non triviale?  A voir quand on écrira la partie de degré 4. 
\end{example} 

\section{Properties of the relative Arason invariant} 
\label{prop.section}
In this section, we study properties of the relative Arason invariant. The results and proofs are largely inspired by~\cite{QT-Arason}, where analogous results were obtained for algebras with orthogonal involutions with trivial discriminant and trivial Clifford invariant. Most of them are contained in~\cite{M:thesis}, under the additional assumption that $\tau_0$ and $\tau$ have split discriminant algebras. 

\subsection{Base point change}
\begin{prop}
\label{Chasles}
Let $\tau_0$, $\tau_1$ and $\tau_2$ be three $F'/F$-unitary involutions on the degree $n$ algebra $B$. We assume all three involutions have isomorphic discriminant algebras if $n$ is even.  
We have 
\[e_3^{\tau_0}(\tau_2)=e_3^{\tau_0}(\tau_1)+e_3^{\tau_1}(\tau_2)\in N^3_{\alpha,\beta}(F),\] 
where $\alpha=0$ if $n$ is odd, $\alpha=[\cd(\tau_0)]=[\cd(\tau_1)]\in\br(F)$ if $n$ is even, and $\beta=[B]\in\br(F')$. 
\end{prop} 
\begin{proof} 
For $i\in\{1,2,3\}$, we let $G_i=\SU(B,\tau_i)$ and $\bar G_i=\PGU(B,\tau_i)$ the corresponding adjoint group. Pick a cocycle $\omega\in Z^1(F,G_0)$ such that its image in $H^1(F,\bar G_0)$ corresponds to the triple $[B,\tau_1,\id_{F'}]$. 
As explained in \S\ref{relative.section}, such a cocycle exists since $G_0$ and $G_1$ have the same Tits class. Moreover, in view of \cite[\S 29]{KMRT}, since any isomorphism $(B,\tau_0)_{F_{\mathrm{sep}}}\rightarrow (B,\tau_1)_{F_{\mathrm{sep}}}$ restricts to an isomorphism $(G_0)_{F_{\mathrm{sep}}}\rightarrow (G_1)_{F_{\mathrm{sep}}}$, the inner twisted form $_{\omega}G_0$ of $G_0$ is isomorphic to $G_1$. We denote by $\phi_\omega$ the induced isomorphism $H^1(F,G_0)\rightarrow H^1(F,G_1)$. 

Now, consider a cocycle $\xi\in H^1(F,G_0)$ with image in $H^1(F,\bar G_0)$ corresponding to $[B,\tau_2,\id_{F'}]$. Again using~\cite[\S 29]{KMRT}, one may check that $\phi_\omega(\xi)$  has image in $H^1(F,\bar G_1)$ also corresponding to $[B,\tau_2,\id_{F'}]$. Moreover, by~\cite[Lemme 7]{Gille}, we have 
\[\rho_{G_1}(\phi_\omega(\xi))=\rho_{G_0}(\xi)-\rho_{G_0}(\omega)\in H^3(F).\] 
Therefore, by definition of the relative Arason invariant, we get 
\[e_3^{\tau_1}(\tau_2)=e_3^{\tau_0}(\tau_2)-e_3^{\tau_0}(\tau_1)\in N^3_{\alpha,\beta}(F),\] 
and this finishes the proof. 
\end{proof} 

\subsection{Relation with the hyperbolic Arason invariant of a sum}
\label{sum.section}
Consider two unitary involutions $\tau_0$ and $\tau$ on a degree $n$ algebra $B$, and pick $s\in \Sym(B,\tau_0)^\times$ such that $\tau=\Int(s)\circ\tau_0$. For any $\lambda\in F^\times$, we let $\theta_\lambda$ be the involution on $M_2(B)$ defined by 
\begin{equation}
\label{sum.eq}
\theta_\lambda\begin{pmatrix} x&y\\z&t\\\end{pmatrix}=\begin{pmatrix} 
\tau_0(x)&-\lambda\tau_0(z)s^{-1}\\
-\lambda^{-1}s\tau_0(y)&s\tau_0(t)s^{-1}\\
\end{pmatrix},\mbox{ for all }x,y,z,t\in B.
\end{equation}
Assume $\tau_0$ is the adjoint involution with respect to some hermitian form $h_0$ with values in the underlying division algebra $D$, endowed with a fixed unitary involution. Then, the involution $\tau$ is adjoint to the hermitian form $h_s$ defined by $h_s(u,v)=h_0(s^{-1}(u),v)$, and $\theta_\lambda$ is adjoint to $h_0\perp\qform{-\lambda}h_s$. 
Therefore, the involution $\theta_{\lambda}$ is called an orthogonal sum of $\tau_0$ and $\tau$. Note that in general, different choices for $\lambda$ may produce non-isomorphic orthogonal sums of $\tau_0$ and $\tau$. Nevertheless, we may use $\theta_\lambda$ to compute $e_3^{\tau_0}(\tau)$ as follows: 
\begin{prop}
\label{sum}
Let $\tau_0$ and $\tau$ be two involutions of the degree $n$ algebra $B$. If $n$ is even, we assume in addition that $\tau_0$ and $\tau$ have isomorphic discriminant algebras. For all $\lambda\in F^\times$, we let $\theta_\lambda$ be the orthogonal sum of $\tau_0$ and $\tau$ defined by~\eqref{sum.eq}. 
\begin{enumerate}
\item If $n$ is odd, then there exists $\lambda\in F^\times$ such that $\theta_\lambda$ has split discriminant algebra. 
If $n$ is even, $\theta_\lambda$ has split discriminant algebra for all $\lambda\in F^\times$. 
\item In both cases, for all $\lambda$ such that $\theta_\lambda$ has split discriminant algebra, we have 
\[e_3^{\tau_0}(\tau)={e_3^{\mathrm{hyp}}(\theta_\lambda)} \in N^3_{\alpha,\beta}(F),\] 
 where   $\alpha=0$ if $n$ is odd, $\alpha=[\cd(\tau_0)]\in\br(F)$ if $n$ is even, and $\beta=[B]=[M_2(B)]\in\br(F')$. 
\end{enumerate} 
\end{prop} 

\begin{remark}
The invariant $e_3^{\mathrm {hyp}}$ has values in $N_{0,\beta}^3(F)$. If $n$ is odd, or $n$ is even and $\cd(\tau_0)$ is split, $e_3^{\tau_0}$ also have values in $N_{0,\beta}^3(F)$, and the proposition states that $e_3^{\tau_0}(\tau)$ and ${e_3^{\mathrm{hyp}}(\theta_\lambda)}$ are equal in this group. Whereas if $n$ is even and $\cd(\tau_0)$ is non-split,  $e_3^{\tau_0}$ has values in $N_{\alpha,\beta}^3(F)$, with $\alpha\not=0$. In this case, the proposition actually states
\[e_3^{\tau_0}(\tau)={e_3^{\mathrm{hyp}}(\theta_\lambda)} \mod F^\times\cdot [\cd(\tau_0)],\]
that is $e_3^{\tau_0}(\tau)$ is the image of $e_3^{\mathrm{hyp}}(\theta_\lambda)$ under 
the projection $N^3_{0,\beta}(F)\rightarrow N^3_{\alpha,\beta}(F)$. 
\end{remark}

\begin{proof} 
Let $\theta_0$ be the involution on $M_2(B)$ defined by 
\[\theta_0\begin{pmatrix} x&y\\z&t\\\end{pmatrix}=\begin{pmatrix} 
\tau_0(x)&-\tau_0(z)\\
-\tau_0(y)&\tau_0(t)
\end{pmatrix},\mbox{ for all }x,y,z,t\in B.\]
It is a hyperbolic involution, since $e=\frac 12 \begin{pmatrix}1&1\\1&1\end{pmatrix}$ is idempotent and satisfies $\theta_0(e)=1-e$, see~\cite[(6.8)]{KMRT}. In particular, its discriminant algebra is split, that is $[\cd(\theta_0)]=0\in\br F$, and $e_3^{\theta_0}=e_3^{\mathrm{hyp}}$. 
Moreover, a direct computation shows that \[\theta_\lambda=\Int\begin{pmatrix}1&0\\0&\lambda s\end{pmatrix}\circ \theta_0.\] 
Therefore, by~\cite[(10.36)]{KMRT}, we have 
\[\cd(\theta_\lambda)\sim\bigl(\delta,\Nrd_{M_2(B)}\begin{pmatrix}1&0\\0&\lambda s\end{pmatrix}\bigr)
\sim \bigl(\delta,\Nrd_B(\lambda s)\bigr).\]
Assume $n=2m+1$ is odd. Since $s$ is symmetric,  $\Nrd_B(s)\in F^\times$. Hence we may chose $\lambda=\Nrd_B(s)$ and we get $\Nrd_B(\lambda s)=\lambda^{n+1}=(\lambda^{m+1})^2$. 
Therefore, $\cd(\theta_\lambda)$ is split for this particular value of $\lambda$. 
Assume now $n=2m$ is even. For all $\lambda\in F^\times$, we have $[\cd(\theta_\lambda)]=\bigr(\delta,\Nrd_B(\lambda s)\bigl)=\bigr(\delta,\lambda^n\Nrd_B(s)\bigl)=\bigr(\delta,\Nrd_B(s)\bigl)$. Since $\tau_0$ and $\tau$ have isomorphic discriminant algebras, we get that $\cd(\theta_\lambda)$ is split, applying again~\cite[(10.36)]{KMRT}. 

This proves the first assertion, and it remains to compare the relative Arason invariants. In the remaining part of the proof, we assume that $\lambda\in F^\times$ has been chosen so that $\cd(\theta_\lambda)$ is split if $n$ is odd. 
Therefore, in both cases, there exists $z\in F'^\times$ such that \[\Nrd_{M_2(B)}\begin{pmatrix}1&0\\0&\lambda s\end{pmatrix}=\Nrd_B(\lambda s)=N_{F'/F}(z).\] It follows that the pair $(\lambda s,z)\in \Sym(B,\tau_0)\times F'^\times$ corresponds to a cohomology class $\xi\in H^1(F,\SU(B,\tau_0))$ under~\eqref{ssym.eq} . Moreover, since $\tau=\Int(s)\circ\tau_0$, the image of $\xi$ in $H^1(F,\PGU(B,\tau_0))$ corresponds to $[B,\tau,\id_{F'}]$. Consider the natural embedding \begin{equation*}i:\,\SU(B,\tau_0)\rightarrow\SU(M_2(B),\theta_0),\ \ \ x\mapsto \begin{pmatrix}1&0\\0&x\end{pmatrix}.\end{equation*} The induced map $i^{(1)}$ in Galois cohomology maps $(\lambda s,z)$ to $(\begin{pmatrix}1&0\\0&\lambda s\end{pmatrix},z)$. Hence, the corresponding element $i^{(1)}(\xi)$ in $H^1\bigl(F,\SU(M_2(B),\theta_0)\bigr)$ maps to the class of $[M_2(B),\theta_\lambda,\id_{F'}]$ in $H^1(F,\PGU(M_2(B),\theta_0))$. Moreover, after scalar extension to an algebraic closure of $F$, the map $i$ corresponds to the standard inclusion $\SL_n\rightarrow \SL_{2n}$, which has Rost multiplier $1$ by~\cite[Example 7.10]{GMS}. 
Therefore, we have 
\[\rho_{SU(M_2(B),\theta_0)}\bigl(i^{(1)}(\xi)\bigr)=\rho_{\SU(B,\tau_0)}(\xi).\]
Assertion (2) follows by definition of the relative Arason invariant. 
\end{proof}

\begin{remark} 
It follows from~\cref{sum} that $e_3^{\mathrm{hyp}}(\theta_{\lambda_1})=e_3^{\mathrm{hyp}}(\theta_{\lambda_2})\in N_{\alpha,\beta}^3(F)$ if $\cd(\theta_{\lambda_i})$ are split for $i=1,2$. If the algebra $B$ is split, this can be checked directly as follows. With the same notations as in the beginning of this section, for $i=1,2$, the involution $\theta_{\lambda_i}$ 
is adjoint to the hermitian form $h_i\simeq h_0\perp\qform{-\lambda_i}h_s$, where $h_0$ and $h_s$ now are hermitian forms of rank $n$ over $(F',\iota)$. If $n$ is odd, $\lambda_1$ and $\lambda_2$ have been chosen so that $h_1$ and $h_2$ have trivial discriminant. If $n$ is even, $h_0$ and $h_s$ have the same discriminant, and so $h_1$ and $h_2$ both have trivial discriminant, for arbitrary values of $\lambda_1$ and $\lambda_2$. In both cases, we have  $e_3^{\mathrm{hyp}}(\theta_{\lambda_i})=e_3(q_{h_i})\in H^3(F)$ for $i=1,2$ by~\cref{splithyp.example}. 
In addition, $h_1-h_2$ is Witt equivalent to $\qform{-\lambda_1,\lambda_2}h_s$, and its trace form $q_{h_1}-q_{h_2}$ is Witt equivalent to 
\[\qform{-\lambda_1,\lambda_2}q_{h_s}\simeq \qform{\lambda_2}\qform{1,-\lambda_1\lambda_2^{-1}}\qform{1,-\delta}\qform{a_1,\dots, a_n},\] where $\qform{a_1,\dots, a_n}$ is any diagonalisation of $h_s$. 

Assume $n=2m+1$ is odd. For $i=1,2$, since $h_i$ has trivial discriminant, $d(h_i)=\lambda_i d(h_0)d(h_s)$ is a norm for the quadratic extension $F'/F$. Hence, 
$(\delta, d(h_1)d(h_2)^{-1})=(\delta, \lambda_1\lambda_2^{-1})=0\in \br(F)$. It follows that the quadratic forms $\qform{1,-\lambda_1\lambda_2^{-1}}\qform{1,-\delta}$ and $q_{h_1}-q_{h_2}$ are hyperbolic. We get 
\[e_3^{\mathrm{hyp}}(\theta_{\lambda_1})-e_3^{\mathrm{hyp}}(\theta_{\lambda_2})=e_3(q_{h_1}-q_{h_2})=0\in H^3(F)\] as expected. 

Assume now $n=2m$ is even, so that $\qform{a_1,\dots, a_n}$ is even dimensional. We have 
\begin{align*}e_3^{\mathrm{hyp}}(\theta_{\lambda_1})-e_3^{\mathrm{hyp}}(\theta_{\lambda_2})=e_3(\qform{-\lambda_1,\lambda_2}q_{h_s})=(\lambda_1\lambda_2^{-1},\delta,(-1)^ma_1\dots a_n).
\end{align*} 
Moreover, by~\cite[(10.35)]{KMRT} \[[\cd(\ad_{h_0})]=(\delta, d(h_0))=(\delta, d(h_s))=(\delta,(-1)^ma_1\dots a_n)\in \br(F).\]
Therefore, given arbitrary $\lambda_1,\lambda_2\in F^\times$, we get  
\[e_3^{\mathrm{hyp}}(\theta_{\lambda_1})-e_3^{\mathrm{hyp}}(\theta_{\lambda_2})=(\lambda_1\lambda_2^{-1})\cdot [\cd(\ad_{h_0})]\in H^3(F).\]
This equality proves $e_3^{\mathrm{hyp}}(\theta_{\lambda_1})$ and $e_3^{\mathrm{hyp}}(\theta_{\lambda_2})$, which generally differ as elements of $H^3(F)$, have the same image in $M^3_\alpha(F)$, as required. 
\end{remark}

\subsection{Unitary involution with a split rank $2$ factor}

Using the same argument as in the proof of~\cref{sum}, we may also compute the hyperbolic Arason invariant of an algebra with unitary involution having a rank $2$ split factor. 

\begin{prop} 
Let $(B,\tau_0)$ be a degree $n$ algebra with unitary involution. 
If $n$ is odd, there exists $\lambda\in F^\times$ such that the unitary involution $\ad_\qform{1,-\lambda}\otimes\tau_0$ has trivial discriminant algebra, and for such a choice of $\lambda$, we have 
\[e_3^{\hyp}(\ad_\qform{1,-\lambda}\otimes\tau_0)=0\in H^3(F)/\cores_{F'/F}(F'^\times\cdot[B]).\]
If $n$ is even, then $\ad_\qform{1,-\lambda}\otimes\tau_0$ has trivial discriminant algebra for all $\lambda\in F^\times$ and 
\[e_3^{\hyp}(\ad_\qform{1,-\lambda}\otimes\tau_0)=(\lambda)\cdot[\cd(\tau_0)]\in H^3(F)/\cores_{F'/F}(F'^\times\cdot[B]).\]
\end{prop} 
\begin{proof}
For all $\lambda\in F^\times$, consider the algebra with involution 
\[(M_2(F'),\ad_{\qform{1,-\lambda}})\otimes_{F'}(B,\tau_0),\] where $\qform{1,-\lambda}$ is a rank $2$ hermitian form over $(F',\iota)$.
It corresponds to $\theta_\lambda$ as defined in the previous section, in the particular case where $\tau=\tau_0$ so that we may assume $s=1$. Hence, \cref{sum} already gives the first part of both assertions, and the formula for $e_3^{\mathrm{hyp}}$ when $n$ is odd. Moreover, since $s=1$, the computation in the proof of~\cref{sum} shows 
that $e_3^{\hyp}(\ad_\qform{1,-\lambda}\otimes\tau_0)$ is given by $\rho_{\SU(B,\tau_0)}(\xi)$, where $\xi\in H^1(F,\SU(B,\tau_0))$ corresponds under~\eqref{ssym.eq} to the class $(\lambda,\lambda^m)$. Therefore, $\xi$ actually belongs to $H^1(F,\mu_{n[F']})$ and the result follows by~\cite[Thm. 1.11]{MPT}, \cite[\S 10]{GM}.  
\end{proof}

\subsection{Order of the relative Arason invariant} 

Using the previous results, we now prove that the relative Arason invariant has order $2$. 
\begin{cor}
\label{order}
Let $\tau_0$ and $\tau$ be two unitary involutions of the degree $n$ algebra $B$. If $n$ is even, we assume in addition that $\tau_0$ and $\tau$ have isomorphic discriminant algebras. We have 
\begin{enumerate}
\item[(a)] $e_3^{\tau_0}(\tau)=e_3^{\tau}(\tau_0)\in N_{\alpha,\beta}^3(F)$, and  
\item[(b)] $2e_3^{\tau_0}(\tau)=0\in N_{\alpha,\beta}^3(F),$
\end{enumerate}
where $\alpha=0$ if $n$ is odd, $\alpha=[\cd(\tau_0)]\in\br F$ if $n$ is even, and $\beta=[B]\in \br(F')$.  
\end{cor}
\begin{proof}
By~\cref{Chasles}, we have $e_3^{\tau_0}(\tau)+e_3^{\tau}(\tau_0)=e_3^{\tau_0}(\tau_0)=0\in N_{\alpha,\beta}^3(F)$. 
Therefore, it is enough to prove assertion (a). This can be checked as follows, using~\cref{sum}. Consider the involution $\theta_\lambda$ defined in~\eqref{sum}. If $n$ is odd, we assume $\lambda$ has been chosen so that $\theta_\lambda$ has split discriminant algebra. As in the beginning of \S\ref{sum.section}, we pick a hermitian form $h_0$ such that $\tau_0\simeq\ad_{h_0}$ and we let $h_s$ be defined by $h_s(u,v)=h_0(s^{-1}(u), v)$, where $s\in \Sym(B,\tau_0)^\times$ is such that $\tau=\Int(s)\circ\tau_0\simeq \ad_{h_s}$. The involution $\theta_\lambda$ is adjoint to $h_0\perp\qform{-\lambda}h_s$. We have 
\[e_3^{\tau_0}(\tau)=e_3^{\mathrm{hyp}}(\theta_\lambda)=e_3^{\mathrm{hyp}}(\ad_{h_0\perp\qform{-\lambda} h_s})\in N_{\alpha,\beta}^3(F).\] On the other hand, since $h_0\perp\qform{-\lambda} h_s$ and $h_s\perp\qform{-\lambda^{-1}}h_0$ are similar hermitian forms, the corresponding involutions are isomorphic and have the same hyperbolic Arason invariant. 
Besides, $\ad_{h_s\perp\qform{-\lambda^{-1}}h_0}$ is an orthogonal sum of $\tau$ and $\tau_0=\Int(s^{-1})\circ\tau$. 
Applying again~\cref{sum}, we thus get 
\[e_3^{\tau_0}(\tau)=e_3^{\mathrm{hyp}}(\ad_{h_s\perp\qform{-\lambda^{-1}}h_0})=e_3^{\tau}(\tau_0)\in N_{\alpha,\beta}^3(F),\] and this concludes the proof. 
\end{proof}

\begin{remark}
As a consequence, we get that given two unitary involutions $\tau_0$ and $\tau$ of $B$, with $\cd(\tau_0)\sim\cd(\tau)$ if $n$ is even, their relative $f_3$ invariant satisfies 
\[f_3^{\tau_0}(\tau)\in \left\{\begin{array}{ll} 
\cores_{F'/F}(F'^\times\cdot [B])&\mbox{ if $n$ is odd,}\\
F^\times\cdot [\cd(\tau_0)]+\cores_{F'/F}(F'^\times\cdot [B])&\mbox{ if $n$ is even.}\\
\end{array}\right.\]
\end{remark}

\section{Comparison with degree $3$ invariants of symplectic involutions} 
\label{symp.section}

The aim of this section is to prove some formulae relating the degree $3$ cohomological invariants in unitary and symplectic type for symplectic descents and quadratic symplectic extensions of a pair of unitary involutions; see Sections~\ref{symplecticdescent.sec} and~\ref{symplecticquadraticextension.sec} below for precise statements. We first recall some known facts on invariants of symplectic involutions. 

\subsection{Degree $3$-invariants for symplectic involutions}
\label{sympinv.sec}

Degree $3$ invariants for symplectic involutions were investigated in~\cite{BMT} and \cite{GPT}. Let $(C,\gamma_0,\gamma)$ be a central simple algebra over $F$ endowed with two symplectic involutions. In particular, $C$ has even degree $r=2\ell$ and exponent dividing $2$. If $\ell$ is even, that is $4\mid r$, a relative invariant $e_3^{\gamma_0}(\gamma)$ is defined in \cite{BMT}. This invariant may also be defined in terms of the Rost invariant of the symplectic group $\Sp(C,\gamma_0)$, see~\cite[\S 5.3]{Tignol:hyderabad}. The process is very similar to the definition of the relative Arason invariant for unitary involutions given in~\cref{relative.section}, with the exact sequence~\eqref{ExactSeq.eq} replaced by
\[1\rightarrow\mu_2\rightarrow\Sp(C,\gamma_0)\rightarrow\PGSp(C,\gamma_0)\rightarrow 1.\] Since we will use it in the proofs of comparison theorems, we briefly recall the procedure. 
The isomorphism class of $(C,\gamma)$ corresponds to a cohomology class $\eta\in H^1\bigr(F,\PGSp(C,\gamma_0)\bigl)$ by~\cite[(29.22)]{KMRT}. Since $\gamma_0$ and $\gamma$ are defined on the same algebra $C$, this class lifts to a class $\xi\in H^1\bigl(F,\Sp(C,\gamma_0)\bigr)$, which is unique up to the action of $H^1(F,\mu_2)$. Moreover, by~\cite[Thm. 1.13]{MPT},\ \cite{GM}, we have 
\[\rho_{\Sp(C,\gamma_0)}(H^1(F,\mu_2))=\left\{\begin{array}{l} 
0 \mbox{ if $\ell$ is even,}\\
F^\times\cdot[C]\mbox{ if $\ell$ is odd.}\\
\end{array}\right.\]
Hence, we may define 
\begin{equation}
e_3^{\gamma_0}(\gamma)=\rho_{\Sp(C,\gamma_0)}(\xi)\in M_\alpha^3(F)=H^3(F)/F^\times\cdot\alpha,
\end{equation}
where $\alpha$ is $0$ if $\ell$ is even and $[C]$ if $\ell$ is odd. 
On the other hand, as explained in~\cite[p. 440]{KMRT}, $\xi$ corresponds to the choice of some $s\in \Sym(C,\gamma_0)^\times$ such that $\gamma=\Int(s)\circ\gamma_0$, and we have 
\[\rho_{\Sp(C,\gamma_0)}(\xi)=\bigl(\nrp_{\gamma_0}(s)\bigr)\cdot[C],\] where 
$\nrp_{\gamma_0}$ denotes the Pfaffian norm map, defined on $\gamma_0$-symmetric elements. 
In particular, we have 
\begin{equation}
\label{RostSymp.eq}
\rho_{\Sp(C,\gamma_0)}(\xi)\in F^\times\cdot[C]\mbox{ for all }\xi\in H^1(F,\Sp(C,\gamma_0)).
\end{equation}
Therefore, $e_3^{\gamma_0}$ is identically zero if $\ell$ is odd. If $\ell$ is even, we recover the original definition given in~\cite{BMT}, namely 
\begin{equation}
\label{RostSympDef.eq}
e_3^{\gamma_0}(\gamma)=\rho_{\Sp(C,\gamma_0)}(\xi)=\bigl(\nrp_{\gamma_0}(s)\bigr)\cdot[C]\in H^3(F),
\end{equation}
where $s\in \Sym(C,\gamma_0)^\times$ satisfies $\gamma=\Int(s)\circ\gamma_0$. 

If the algebra $C$ has even co-index, it is endowed with a hyperbolic involution, and the construction above leads to a hyperbolic invariant $e_3^{\hyp}(\gamma)\in H^3(F)$, which may be considered as an absolute invariant of $\gamma$. Unexpectedly, this invariant extends to a more general setting. In~\cite{GPT}, Garibaldi, Parimala and Tignol proved there exists a degree $3$ invariant $e_3(\gamma)\in H^3(F)$ for all symplectic involution $\gamma$ on an algebra $C$ of degree divisible by $8$, including division algebras. This invariant vanishes on hyperbolic involutions, and it satisfies $e_3^{\gamma_0}(\gamma)=e_3(\gamma)-e_3(\gamma_0)$ for all symplectic involutions $\gamma_0$ and $\gamma$ of $C$. 

\begin{example}
\label{SympWithOrthDescent}
Assume there exists an even degree central simple algebra $A$ over $F$, endowed with two orthogonal involutions $\sigma_0$ and $\sigma$, and a quaternion algebra $Q$ with canonical involution $\ba$ such that 
\[(C,\gamma_0)\simeq (A,\sigma_0)\otimes_F (Q,\ba)\mbox{ and }(C,\gamma)\simeq (A,\sigma)\otimes_F (Q,\ba).\]
Pick $s\in \Sym(A,\sigma_0)^\times$ such that $\sigma=\Int(s)\circ\sigma_0$; then $s\otimes 1\in \Sym(C,\gamma_0)^\times$ and $\gamma=\Int(s\otimes 1)\circ\gamma_0$. Moreover, by~\cite[Lemma 9(d)]{BMT}, and~\cite[(7.3)(1)]{KMRT}
\[\nrp_{\gamma_0}(s\otimes 1)=\Nrd_A(s)=d(\sigma)d(\sigma_0)^{-1}\in F^\times/F^{\times 2}.\]
Applying the projection formula and~\cite[Cor. 2.6.10]{GS}, we have $\Nrd_A(s)\cdot[A]=0$. Therefore, since $[C]=[A]+[Q]$, we get 
\[e_3^{\gamma_0}(\gamma)=\bigl(d(\sigma)d(\sigma_0)^{-1}\bigr)\cdot[Q]\in H^3(F).\]
\end{example}

\subsection{Unitary involutions with a symplectic descent}
\label{symplecticdescent.sec} 
In this section, we consider a central simple algebra $B$ with two $F'/F$-unitary involutions $\tau_0$ and $\tau$, and we assume they both admit a symplectic descent to the same algebra $C$ over $F$, that is 
\begin{equation}
\label{sympdescent.eq}
(B,\tau_0)\simeq (C,\gamma_0)\otimes(F',\iota)\mbox{ and }(B,\tau)=(C,\gamma)\otimes(F',\iota),
\end{equation}
for some symplectic involutions $\gamma_0$ and $\gamma$ of $C$. 

\begin{remark}
\label{existencedescent}
 Condition~(\ref{sympdescent.eq}) is a strong condition on the triple $(B,\tau_0,\tau)$, which does not hold in general.  
First of all, the existence of a symplectic descent implies that $B$ has even degree $n=2m$, and exponent dividing $2$. Conversely, if the algebra $B$ satisfies those conditions, it admits a descent up to Brauer equivalence since it has trivial corestriction, but need not have a descent up to isomorphism, see~\cite[Rem. 4.8]{Bar14} and~\cref{hypnotdec}. 
In addition, even assuming $B\simeq C\otimes_F F'$ for some central simple algebra $C$ over $F$, a unitary involution $\tau$ of $B$ generally does not have a symplectic descent to $C$. 
Assume for instance $B$ is a quaternion algebra $Q$ with trivial corestriction, and pick a particular quaternion descent $Q_0$ over $F$. The algebra $Q_0$ is not unique in general, and a unitary involution $\tau$ on $Q$ has a symplectic descent to $Q_0$ if and only if $Q_0$ is the discriminant algebra of $\tau$, see~\cite[(2.22)\& p. 129]{KMRT}. By uniqueness of the symplectic involution on a quaternion algebra, our condition holds in degree $2$ if and only if the unitary involutions $\tau_0$ and $\tau$ are isomorphic. 
\end{remark}

Under condition~(\ref{sympdescent.eq}), the invariants of the pairs $(\gamma_0,\gamma)$ and $(\tau_0,\tau)$ are related as follows :
\begin{prop}
\label{symplecticdescent.prop}
Let $(B,\tau_0,\tau)$ be a degree $n$ central simple algebra with two $F'/F$-unitary involutions. We assume there exists a central simple algebra over $F$ with two symplectic involutions $(C,\gamma_0,\gamma)$, such that 
\[(B,\tau_0)\simeq (C,\gamma_0)\otimes(F',\iota)\mbox{ and }(B,\tau)=(C,\gamma)\otimes(F',\iota).\]
Then, $n=2m$ is even, $B$ has exponent dividing $2$, $\cd(\tau_0)\sim\cd(\tau)\sim C^{\otimes m}$, and 
\[e_3^{\tau_0}(\tau)=
\left\{
\begin{array}{ll} 
e_3^{\gamma_0}(\gamma)\in H^3(F)/N_{F'/F}(F'^\times)\cdot[C]&\mbox{ if $m$ is even,}\\
0\in H^3(F)/F^\times\cdot[C]&\mbox{ if $m$ is odd.}\\
\end{array}\right.\]
Moreover, in both cases, we have $f_3^{\tau_0}(\tau)=0$. 
\end{prop}

\begin{remark}
If $m$ is even, the invariant $e_3^{\gamma_0}(\gamma)$ has values in $H^3(F)$, and the proposition says its image in the quotient $H^3(F)/N_{F'/F}(F'^\times)\cdot[C]$ coincides with $e_3^{\tau_0}(\tau)$. 
\end{remark}

\begin{proof}
Since $C$ admits symplectic involutions, the algebras $C$ and $B$ have even degree $n=2m$, and exponent dividing $2$. 

To compute the Brauer class of $\cd(\tau_0)$, we may assume the algebra $B$ is split, by the argument given in~\cite[Proof of (10.36)]{KMRT}. 
Under this assumption, $C$ is split by $F'$, hence it is Brauer equivalent to a quaternion algebra $Q=(\delta, a)$ for some $a\in F^\times$. It follows that $(C,\gamma_0)$ decomposes as $(M_m(F),\ad_{\varphi_0})\otimes(Q,\ba)$ for some $m$-dimensional quadratic form $\varphi_0$ over $F$, see~\cite[\S5.1]{Tignol:hyderabad} or~\cite[Prop. 2.1(2)]{QT-conic}, and where $\ba$ denotes the canonical involution of $Q$.  
Therefore, by~\eqref{QuatUnit.eq}, we have 
\[(B,\tau_0)\simeq (M_m(F),\ad_{\varphi_0})\otimes (Q,\ba)\otimes (F',\iota)\simeq (M_n(F),\ad_{\varphi_0\otimes\pform{a}})\otimes(F',\iota).\]
So~\eqref{DiscOrthDescent.eq} gives $\cd(\tau_0)\sim(\delta, a^m)\sim C^{\otimes m}$, as required. 
The same argument also applies to the involution $\tau$, and it remains to compare the relative $e_3$ invariants. 

Note first that the invariant $e_3^{\tau_0}$ has values in $N_{\alpha,\beta}^3(F)$, with $\alpha=0$ if $m$ is even, and $\alpha=[C]$ if $m$ is odd, and $\beta=[B]$.  
Morever, by the projection formula, for all $z\in F'^\times$, we have $\cores_{F'/F}\bigl((z)\cdot\beta\bigr)=N_{F'/F}(z)\cdot [C]$, which clearly belongs to $F^\times \cdot [C]$. Therefore, $e_3^{\tau_0}$ has values in $H^3(F)/N_{F'/F}(F'^\times)\cdot[C]$ if $m$ is even and 
$H^3(F)/F^\times\cdot[C]$ if $m$ is odd. 

To conclude the proof, we compare the Rost invariants. Since \[(B,\tau_0)\simeq (C,\gamma_0)\otimes (F',\iota),\] there is a natural embedding $i:\,\Sp(C,\gamma_0)\rightarrow \SU(B,\tau_0)$. When both groups are split, identifying $\SU(M_n(F)\times M_n(F)^{\mathrm{op}},\varepsilon)$ with $\SL_n(F)$, this map corresponds to the canonical embedding $\Sp_n\rightarrow \SL_n$, which has Rost multiplier $1$ by~\cite[Example 7.11]{GMS}. Therefore, for all $\xi\in H^1(F,\Sp(C,\gamma_0))$, we have \[\rho_{\SU(B,\tau_0)}(i^{(1)}(\xi))=\rho_{\Sp(C,\gamma_0)}(\xi).\]

Consider $s\in\Sym(C,\gamma_0)^\times$ such that $\gamma=\Int(s)\circ \gamma_0$. 
It determines a class $\xi\in H^1\bigl(F,\Sp(C,\gamma_0)\bigr)$ with image in $H^1\bigl(F,\PGSp(C,\gamma_0)\bigr)$ corresponding to the isomorphism class of $(C,\gamma)$. Moreover, since $\tau=\Int(s\otimes 1)\circ\tau_0$, the element $i^{(1)}(\xi)\in H^1\bigl(F,\SU(B,\tau_0)\bigr)$ has image in $H^1\bigl(F,\PGU(B,\tau_0)\bigr)$ corresponding to $[B,\tau,\id_{F'}]$. 
Therefore,
\[e_3^{\tau_0}(\tau)=\rho_{\SU(B,\tau_0)}(i^{(1)}(\xi))=\rho_{\Sp(C,\gamma_0)}(\xi)\in N_{\alpha,\beta}^3(F).\]
The main result follows by~\eqref{RostSymp.eq} if $m$ is odd and~\eqref{RostSympDef.eq} if $m$ is even. Since $e_3^{\gamma_0}(\gamma)$ has values in $_2H^3(F)$, we also get $f_3^{\tau_0}(\tau)=0$. 
\end{proof}

\begin{example}
Assume $B$ is split; then $C$ is Brauer equivalent to a quaternion algebra $Q=(\delta,a)$ for some $a\in F^\times$. 
Applying as above~\cite[\S5.1]{Tignol:hyderabad} or~\cite[Prop. 2.1(2)]{QT-conic}, we get 
\[(C,\gamma_0)=(M_m(F),\ad_{\varphi_0})\otimes (Q,\ba)\mbox{ and }(C,\gamma)=(M_m(F),\ad_{\varphi})\otimes (Q,\ba).\]
Hence $e_3^{\gamma_0}(\gamma)=(d(\varphi)d(\varphi_0)^{-1})\cdot [Q]$ by~\cref{SympWithOrthDescent}. 
Moreover, since $B$ is split, $N_{F'/F}(F'^\times)\cdot [C]=\cores_{F'/F}(F'^\times\cdot[B])=0$. Therefore, we have \[e_3^{\tau_0}(\tau)=
\left\{
\begin{array}{ll}
\bigr(d(\varphi)d(\varphi_0)^{-1}\bigl)\cdot[Q]\in H^3(F)&\mbox{ if $m$ is even,}\\
0\in H^3(F)/F^\times\cdot[Q]&\mbox{ if $m$ is odd.}\\
\end{array}\right.\]
Of course, this is compatible with the result given by~\cref{split}. 
This can be checked using~\eqref{QuatUnit.eq}, which shows that the traces of the hermitian forms associated to $\tau_0$ and $\tau$ are respectively similar to $\qform{1,-\delta}\qform{1,-a}\varphi_0$ and $\qform{1,-\delta}\qform{1,-a}\varphi$. 
\end{example}
\subsection{Symplectic quadratic extensions of unitary involutions}
\label{symplecticquadraticextension.sec}
The algebra with symplectic involution $(C,\gamma)$ is called a quadratic extension of the algebra with unitary involution $(B,\tau)$ if $C$ contains a $\gamma$-stable quadratic \'etale $F$-algebra $\tilde F$ such that 
\[(\tilde F,\gamma_{\mid \tilde F})\simeq (F',\iota)\mbox{ and }\bigr(Z_C(\tilde F),\gamma_{\mid Z_C(\tilde F)}\bigl)\simeq (B,\tau),\] where $Z_C(\tilde F)$ denotes the centraliser of $\tilde F$ in $C$. When this holds, $C$ is a central simple $F$-algebra of degree $2n$, where $n$ is the degree of $B$, and is satisfies $C_{F'}\sim B$. In particular, this implies $B$ has exponent at most $2$. 

By the Skolem-Noether theorem, if $C$ contains two distinct subalgebras $\tilde F_1$ and $\tilde F_2$ as above, the corresponding centralisers $Z_C(\tilde F_1)$ and $Z_C(\tilde F_2)$ are isomorphic. Nevertheless, the restrictions of $\gamma$ to both centralisers are not isomorphic in general. Therefore, a given $(C,\gamma)$ might be a quadratic extension of some $(B,\tau_1)$ and $(B,\tau_2)$ for non isomorphic unitary involutions $\tau_1$ and $\tau_2$. For instance, the (hyperbolic) algebra with symplectic involution $(M_4(F),\gamma)$ is a quadratic extension of 
$(M_2(F),\ba)\otimes_F (F',\iota)$ and of $\bigl((\delta, a),\ba\bigr)\otimes_F(F',\iota)$ by~\cite[Example 1.13]{GaQ:deg12}. 
In both cases, the algebra is $M_2(F')\simeq (\delta,a)\otimes_F F'$, nevertheless, the involutions are non isomorphic if $(\delta, a)$ is non split, since they have distinct discriminant algebras. 

As opposed to this, the involution $\gamma$ is uniquely determined by its restriction $\tau$. More precisely, given a degree $n$ central simple algebra with $F'/F$-unitary involution $(B,\tau)$ such that $B$ has exponent at most $2$, and a central simple algebra $C$ of degree $2n$ over $F$ such that $C_{F'}\sim B$, there exists a unique symplectic involution $\gamma$ of $C$ such that $(C,\gamma)$ is a quadratic extension of $(B,\tau)$. This is an easy consequence of the Skolem-Noether theorem, see~\cite[Prop. 1.9]{GaQ:deg12}. We call $\gamma$ the symplectic extension of $\tau$ to $C$. 

\begin{prop}
\label{symplecticextension}
Let $(B,\tau_0,\tau)$ be a degree $n$ algebra with two $F'/F$ unitary involutions. We assume $B$ has exponent dividing $2$, and if $n$ is even, we assume $\tau_0$ and $\tau$ have isomorphic discriminant algebras. Let $C$ be a degree $2n$ algebra over $F$ such that $C_{F'}\sim B$, and denote by $\gamma_0$, respectively $\gamma$, the symplectic extension of $\tau_0$, respectively $\tau$, to $C$. Their relative Arason invariant is given by
\[e_3^{\gamma_0}(\gamma)=f_3^{\tau_0}(\tau)\in H^3(F).\]
\end{prop} 

\begin{proof} 
The proof is similar to the proof of~\cref{symplecticdescent.prop}, with the role of the symplectic and unitary involutions reversed. 
Since $(C,\gamma_0)$ is a quadratic symplectic extension of $(B,\tau_0)$, there is a canonical embedding $i:\,SU(B,\tau_0)\rightarrow \Sp(C,\gamma_0)$. 
Pick a cohomology class $\xi\in H^1\bigl(F,\SU(B,\tau_0)\bigr)$ with image in $H^1\bigl(F,\PGU(B,\tau_0)\bigr)$ corresponding to $[B,\tau,\id_{F'}]$. 
Such a class exists since $\tau_0$ and $\tau$ have isomorphic discriminant algebras if $n$ is even. 
By~\eqref{ssym.eq}, $\xi$ corresponds to the class of a pair $(s,z)\in \SSym(B,\tau_0)$, where $s\in\Sym(B,\tau_0)^\times$ satisfies $\tau=\Int(s)\circ\tau_0$. Clearly, $\Sym(B,\tau_0)\subset\Sym(C,\gamma_0)$. 
Therefore, $\Int(s)\circ\gamma_0$ is a symplectic involution on $C$, acting as $\tau$ on $B$. By uniqueness of the quadratic symplectic extension, we thus have $\gamma=\Int(s)\circ \gamma_0$. 
Hence, we get 
\[f_3^{\tau_0}(\tau)=2\rho_{\SU(B,\tau_0)}(\xi)\in H^3(F)\mbox{ and }e_3^{\gamma_0}(\gamma)=\rho_{\Sp(C,\gamma_0)}\bigl(i^{(1)}(\xi)\bigr)\in H^3(F).\] Therefore, the following lemma finishes the proof : 
\begin{lem}
The natural embedding $i:\SU(B,\tau_0)\rightarrow \Sp(C,\gamma_0)$ has Rost multiplier $2$, that is for all $\xi\in H^1\bigl(F,\SU(B,\tau_0)\bigr)$, we have 
$\rho_{\Sp(C,\gamma_0)}(i^{(1)}(\xi))=2\rho_{\SU(B,\tau_0)}(\xi)$. 
\end{lem}
To compute the Rost multiplier, we use the same strategy as for the examples given in~\cite[\S 7]{GMS}. First of all, by~\cite[Prop. 7.9(4)]{GMS}, we may assume both groups are split. 
In other words, viewing matrices in $M_{2n}(F)$ as $2\times 2$ block matrices, we may identify $(C,\gamma_0)$ with 
$M_{2n}(F)$ endowed with the involution defined by  
\[\gamma_0\begin{pmatrix}A&B\\C&D\end{pmatrix}=\begin{pmatrix}D^t&-B^t\\-C^t&A^t\end{pmatrix},\mbox{ for all }A,B,C,D\in M_n(F).\]
Therefore, the map 
\[\psi:\,\bigl(M_n(F)\times M_n(F)^{\mathrm{op}},\varepsilon)\rightarrow (M_{2n}(F),\gamma_0),\ \ 
(X,Y^{\mathrm{op}})\mapsto \begin{pmatrix}X&0\\0&Y^t\end{pmatrix}\]
preserves the involutions. Moreover, $\tilde F=F[\Delta]$, where $\Delta=\begin{pmatrix}1&0\\0&-1\end{pmatrix}$, satisfies 
$(\tilde F,\gamma_{\mid\tilde F})\simeq (F\times F, \iota)$, and the centraliser, endowed with the restriction of $\gamma_0$, coincides with the image of $\psi$. Hence, $\gamma_0$ is the symplectic quadratic extension of $\varepsilon$ to $M_{2n}(F)$ and the map $i$ we are interested in corresponds in the split case to the restriction of $\psi$ to the group 
$\SU\bigl(M_n(F)\times M_n(F)^{\mathrm{op}},\varepsilon\bigr)$. 
On the other hand, this group is isomorphic to $\SL_n(F)$ by $X\mapsto (X,(X^{-1})^{\mathrm {op}})$. 
So we are computing the Rost multiplier of the embedding 
\[\SL_n(F)\rightarrow \Sp_{2n}(F), \ \ X\mapsto \begin{pmatrix} X&0\\0&(X^{-1})^{t}\end{pmatrix}.\]

By~\cite[Example 7.11]{GMS}, the standard embedding $\Sp_{2n}\rightarrow \SL_{2n}$ has Rost multiplier $1$, therefore, we may replace $\Sp_{2n}$ by $\SL_{2n}$ by~\cite[Prop. 7.9(1)]{GMS}. The result now follows from~\cref{rostmultiplier}. 

\end{proof} 

\section{Comparison with degree $3$ invariants of orthogonal involutions} 
\label{orth.section}

We now prove some formulae relating the degree $3$ cohomological invariants of orthogonal descents and quadratic orthogonal extensions of pairs of unitary involutions; see Sections~\ref{orthogonalquadraticextension.sec} and~\ref{orthogonaldescent.sec} for precise statements. First, we extend the definition of the relative Arason invariant for orthogonal involutions to involutions with isomorphic non split Clifford algebras, see Section \ref{orthogonalinvariant.sec} below. The definition of this relative invariant follows the same lines as in the unitary and symplectic cases, but there are some complications if the algebra has even degree, related to the fact that given two involutions with the same discriminant, there are two different ways of identifying the centers of their Clifford algebras. Therefore, an additional condition is required for the relative Arason invariant to be defined, see~\eqref{addcond.eq}. This condition can also be interpreted in terms of outer automorphisms of groups of type ${\mathrm D}$, following~\cite{QT-AutExt}, and in terms of similarity factors of hermitian forms, see Remarks~\ref{relorth.rem}(2) and~\ref{sim.rem}. 

\subsection{Degree $3$-invariants for orthogonal involutions}
\label{orthogonalinvariant.sec}

Using the Rost invariant, one may define absolute and relative Arason invariants for some orthogonal involutions, see~\cite{Tignol:hyderabad} and~\cite{QT-Arason} for precise definitions and some properties. The purpose of this subsection is to slightly expand this definition, by considering pairs of involutions having isomorphic (and generally non trivial) Clifford algebras. In particular, if the underlying algebra has even degree, this implies that the involutions have the same discriminant. 
Throughout, we assume $\sigma_0$ and $\sigma$ are two orthogonal involutions of the central simple algebra $A$. 

\subsubsection{Odd degree} 
\label{odddegree.sec}
Assume first that $A$ has odd degree $r=2\ell+1$, with $\ell\geq 2$, so that the underlying algebraic groups have type ${\mathrm B}_\ell$. Then the algebra $A$ is split, and there exists two quadratic forms $\varphi_0$ and $\varphi$, unique up to a scalar factor, such that $\sigma_0=\ad_{\varphi_0}$ and $\sigma=\ad_{\varphi}$. 
We assume in addition that $\varphi_0$ and $\varphi$ have isomorphic even Clifford algebras. 
Under this condition, we have 
\[\varphi_0-\qform{\disc(\varphi_0)\disc(\varphi)^{-1}}\varphi\in I^3(F),\]
and we may define a relative Arason invariant by
\begin{equation}
\label{OrthOddDeg}
e_3^{\sigma_0}(\sigma)=e_3\bigl(\varphi_0-\qform{\disc(\varphi_0)\disc(\varphi)^{-1}}\varphi\bigr)\in\,_2H^3(F).
\end{equation}
Indeed, the quadratic form $\varphi_0-\qform{\disc(\varphi_0)\disc(\varphi)^{-1}}\varphi$ has dimension $2r$ and trivial discriminant. Moreover, by~\cite[Chap. V (3.13)\&(3.16)]{Lam}, its Clifford invariant is $[\C_0(\varphi_0)]+[\C_0(\varphi)]\in\br(F)$, which vanishes if and only if $\C_0(\varphi_0)$ and $\C_0(\varphi)$ are isomorphic, since those algebras have the same degree and exponent at most $2$. 
Therefore, we may consider the Arason invariant of this quadratic form
\[e_3\bigl(\varphi_0-\qform{\disc(\varphi_0)\disc(\varphi)^{-1}}\varphi\bigr)\in H^3(F).\]
The value of this invariant does not depend on the particular choices we made for $\varphi_0$ and $\varphi$ in their respective similarity classes. Indeed, for all $\lambda_0,\lambda\in F^\times$, we have 
\[\qform{\lambda_0}\varphi_0-\qform{\disc(\qform{\lambda_0}\varphi_0)\disc(\qform{\lambda}\varphi)^{-1}}\qform{\lambda}\varphi=\qform{\lambda_0}(\varphi_0-\qform{\disc(\varphi_0)\disc(\varphi)^{-1}}\varphi),\] and similar forms in $I^3(F)$ have the same Arason invariant. 

\begin{remark}
(1) Alternately, one may chose $\varphi_0$ and $\varphi$ of discriminant $1$, in which case they are uniquely defined up to isometry, and we have \[e_3^{\sigma_0}(\sigma)=e_3(\varphi_0-\varphi)\in H^3(F).\]

(2) In view of the description of the Rost invariant for the group $\Spin(\varphi_0)$ given in~\cite[p. 437]{KMRT}, we also have 
\[e_3^{\sigma_0}(\sigma)=\rho_{\Spin(\varphi_0)}(\xi),\] where $\xi\in H^1\bigl(F,\Spin(\varphi_0)\bigr)$ is any cocycle with image in $H^1(F,\Orth(\varphi_0))$ corresponding under~\cite[(29.28)]{KMRT} to the class of a quadratic form $\varphi$ such that $\sigma\simeq\ad_\varphi$. Rephrasing this in terms of involutions, $\xi$ is any cocycle with image in $H^1\bigr(F,\Orth(A,\sigma_0)\bigl)$ corresponding under~\cite[(29.26)]{KMRT} to the class of a symmetric element $s\in \Sym(A,\sigma_0)$ such that $\sigma=\Int(s)\circ\sigma_0$, or equivalently any cocycle with image in $H^1\bigr(F,\GO(A,\sigma_0)\bigl)$ corresponding under the bijection described in~\cite[p. 405]{KMRT} to the conjugacy class of $\sigma$. 
\end{remark}

\subsubsection{Even degree} Assume now $A$ has even degree $r=2\ell$, with $\ell\geq 3$, so that the underlying groups have type ${\mathrm D}_\ell$. Throughout this subsection, we assume in addition that the Clifford algebras $\C(A,\sigma_0)$ and $\C(A,\sigma)$ are $F$-isomorphic. In particular, they have isomorphic center $Z_0\simeq Z$, so that $\sigma_0$ and $\sigma$ have the same discriminant. 

Pick an $F$-isomorphism 
\[\Psi:\,\C(A,\sigma) \rightarrow \C(A,\sigma_0),\] and denote by $\psi:\,Z\rightarrow Z_0$ the induced isomorphism of the centers. The triple $(A,\sigma,\psi)$ corresponds to an element $\eta\in H^1\bigl(F,\PGO^+(A,\sigma_0)\bigr)$ by~\cite[\S 29.F]{KMRT}. Let $\mu_0$ and $\mu$ denote the respective centers of the groups $\Spin(A,\sigma_0)$ and $\Spin(A,\sigma)$. We have 
\[\mu_0=\left\{\begin{array}{ll} 
\mu_{4[Z_0]}&\mbox{ if $\ell$ is odd,}\\
R_{Z_0/F}(\mu_2)&\mbox{ if $\ell$ is even,}\\
\end{array}\right.\]
and similarly for $\mu$. Therefore, the isomorphism $\psi$ induces an isomorphism between $\mu$ and $\mu_0$, hence also an isomorphism $\psi^{(2)}:\,H^2(F,\mu)\rightarrow H^2(F,\mu_0)$. We may use $\psi^{(2)}$ to compare the Tits classes of the groups $\Spin(A,\sigma)$ and $\Spin(A,\sigma_0)$. 

If $\ell$ is even, so that $H^2(F,\mu_0)\simeq \,_2\br(Z_0)$, we have $t_{\Spin(A,\sigma_0)}=[\C(A,\sigma_0)]$ by ~\cite[(31.13)]{KMRT}. Assume now $\ell$ is odd. We may identify $H^2(F,\mu_0)$ with a subgroup of $\br(F)\times\br(Z_0)$ by ~\cite[Prop. 2.10]{CTGP}, and it follows from~\cite[(31.11)]{KMRT} that the Tits class of $\Spin(A,\sigma_0)$ corresponds to the pair $\bigl([A],[\C(A,\sigma_0)]\bigr)$. In both cases, the $F$-linear isomorphism $\Psi$ induces a $Z_0$-linear isomorphism between 
\[\C(A,\sigma)\otimes_{Z,\psi} Z_0\rightarrow \C(A,\sigma_0),\] where the tensor product structure is given by $\psi:\,Z\rightarrow Z_0$. Therefore, we have \[\psi^{(2)}\bigr([\C(A,\sigma)]\bigl)=[\C(A,\sigma)\otimes_{Z,\psi} Z_0]=[\C(A,\sigma_0)]\in\br(Z_0).\]
Hence, independently of the parity of $\ell$, the Tits classes of $\Spin(A,\sigma_0)$ and $\Spin(A,\sigma)$ coincide under the identification $\psi^{(2)}$. 

With this in hand, we aim at defining a relative Arason invariant $e_3^{\sigma_0}(\sigma)$ as we did in the unitary case, replacing~\eqref{ExactSeq.eq} by 
\[1\rightarrow \mu_0\rightarrow \Spin(A,\sigma_0)\rightarrow \PGO^+(A,\sigma_0)\rightarrow 1.\]
This sequence induces a connecting map 
\[\partial:\,H^1\bigr(F,\PGO^+(A,\sigma_0)\bigl)\rightarrow H^2(F,\mu_0).\]
Applying again a twisting argument as in~\cite[I.5.4]{Serre} (see also~\cite[\S 1]{Gar}), one may check that $\partial(\eta)$ is the difference between the Tits classes of the groups $\Spin(A,\sigma)$ and $\Spin(A,\sigma_0)$, where the first Tits class, which belongs to $H^2(F,\mu)$ is viewed as an element of $H^2(F,\mu_0)$ via the isomorphism $\psi^{(2)}$. 
Hence, as we have just checked, we have $\partial(\eta)=0$, so that $\eta$ lifts to a class $\xi\in H^1\bigr(F,\Spin(A,\sigma_0)\bigl)$, which is unique up to the action of $H^1(F,\mu_0)$. 
In view of the behavior of the Rost invariant under twisting as described in~\cite[Lemme 7]{Gille} (see also~\cite[Prop. 1.7]{MPT}), it follows that the class $\eta$ determines $\rho_{\Spin(A,\sigma_0)}(\xi)\in H^3(F)$ up to an element which belongs to the image of $H^1(F,\mu_0)$ under $\rho_{\Spin(A,\sigma_0)}$. This image is described in~\cite[Thm. 1.15]{MPT},\ \cite{GM}, and we have 
\[\rho_{\Spin(A,\sigma_0)}\bigr(H^1(F,\mu_0)\bigl)\subset
\left\{\begin{array}{ll}
F^\times\cdot [A]+\cores_{Z_0/F}\bigl(Z_0^\times\cdot [\C(A,\sigma_0)]\bigr)&\mbox{ if $\ell$ is odd,}\\
\cores_{Z_0/F}\bigl(Z_0^\times\cdot [\C(A,\sigma_0)]\bigr)&\mbox{ if $\ell$ is even.}\\
\end{array}\right.\]
Note that if $\ell$ is even, then by~\cite[(9.14)]{KMRT}, we have $\cores_{Z_0/F}([\C(A,\sigma_0)])=[A].$ Therefore, by the projection formula, for all $\lambda\in F^\times$, we have  \[(\lambda)\cdot[A]=\cores_{Z_0/F}\bigl((\lambda)\cdot[\C(A,\sigma_0)]\bigl).\] Hence, independently of the parity of $\ell$, we actually have 
\begin{equation}
\label{RostCenterOrth.eq}
\rho_{\Spin(A,\sigma_0)}\bigr(H^1(F,\mu_0)\bigl)\subset F^\times\cdot [A]+\cores_{Z_0/F}\bigl(Z_0^\times\cdot [\C(A,\sigma_0)]\bigr). 
\end{equation}

We still have to check that the value of $\rho_{\Spin(A,\sigma_0)}(\xi)$ in the quotient of $H^3(F)$ by this subgroup depends only on the involution $\sigma$ of $A$, and not of the cocycle $\eta$, which corresponds to the triple $(A,\sigma,\psi)$. An additional condition is required to guarantee this fact. More precisely, assume there exists an $F$-isomorphism \[\Psi':\,\C(A,\sigma) \rightarrow \C(A,\sigma_0),\] such that the induced isomorphism on the centers is $\psi\circ\iota:\, Z\rightarrow Z_0$, where $\iota$ denotes here the non-trivial $F$-automorphism of $Z$. Composing $\Psi'$ with $\Psi^{-1}$, we get that the algebra $\C(A,\sigma)$ is isomorphic to its conjugate $^{\iota}\C(A,\sigma)$. Under this condition, the class $\eta'\in H^1\bigl(F,\PGO^+(A,\sigma_0)\bigr)$ corresponding to the triple $(A,\sigma,\iota\circ\psi)$ also lifts to a $\xi'\in H^1\bigr(F,\Spin(A,\sigma_0)\bigl)$, and it is not known in general how the values of the Rost invariant compare for $\xi$ and $\xi'$. 
Hence we have to stay away from this situation. 
So, from now on, we assume 
\begin{multline}
\label{addcond.eq}
\mbox{Either $\C(A,\sigma)$ and its conjugate $^\iota\C(A,\sigma)$ are not isomorphic, }\\
\mbox{or $(A,\sigma)$ admits improper similitudes.}
\end{multline} 
In the first case, the class $\eta'$ satisfies $\partial(\eta')\not =0$, so $\eta'$ does not lift to an element $\xi'\in H^1\bigl(F,\Spin(A,\sigma_0)\bigr)$. In the second case, any improper similitude $g$ for $(A,\sigma)$ induces an automorphism $\Int(g)$ of $(A,\sigma)$ and an automorphism $\C(g)$ of $\C(A,\sigma)$ which acts non trivially on $Z$ by~\cite[(13.2)]{KMRT}. Therefore, the triples $(A,\sigma,\psi)$ and $(A,\sigma,\iota\circ \psi)$ are $F$-isomorphic, and it follows that $\eta=\eta'\in H^1\bigr(F,\PGO^+(A,\sigma_0)\bigl)$. Hence, in both cases, there is a unique $\eta\in H^1\bigr(F,\PGO^+(A,\sigma_0)\bigl)$ which is associated to $\sigma$ and which lifts to a $\xi\in H^1\bigl(F,\Spin(A,\sigma_0)\bigr)$. 
So we have a well defined invariant 
\[e_3^{\sigma_0}(\sigma)=\rho_{\Spin(A,\sigma_0)}(\xi)\in N^3_{\alpha,\beta}(F),\] 
where $\alpha=[A]$, $\beta=[\C(A,\sigma_0)]\in \br(Z_0)$, and the element $\xi\in H^1\bigr(F,\Spin(A,\sigma_0)\bigl)$ has image in $H^1(F,\PGO^+\bigl(A,\sigma_0)\bigr)$ corresponding to the triple $(A,\sigma,\psi)$. 
We will refer to $e_3^{\sigma_0}$ as the relative Arason invariant with $\sigma_0$ as a base point. It has values in 
\[N^3_{\alpha,\beta}(F)=H^3(F)/\bigl(F^\times\cdot[A]+\cores_{Z_0/F}(Z_0^\times\cdot[\C(A,\sigma_0)])\bigl),\]
which coincides with $N^3_{0,\beta}(F)=H^3(F)/\cores_{Z_0/F}(Z_0^\times\cdot[\C(A,\sigma_0)])$ when $\ell$ is even. 

\begin{remark}
\label{relorth.rem}
\begin{enumerate}
\item Using the so-called fundamental relations~\cite[(9.12)]{KMRT}, one may check that $\C(A,\sigma)$ is isomorphic to its conjugate $^\iota\C(A,\sigma)$ if and only if the algebra $A$ is split by the quadratic extension $Z$ corresponding to the discriminant of $\sigma$. 
\item Condition~\eqref{addcond.eq} can be rephrased in the language of~\cite{QT-AutExt}; it is equivalent to saying that either the group $G=\PGO^+(A,\sigma)$ admits an outer automorphism defined over $F$, or its Tits class is not fixed under the nontrivial automorphism of its Dynkin diagram. We need to avoid groups for which there is no Tits class obstruction, and yet, there is no outer automorphism defined over $F$, see loc. cit. Prop. 2.5 and Thm 1.1(2). 
\item Assume $\sigma_0$ and $\sigma$ have trivial discriminant, so that $Z\simeq F\times F\simeq Z_0$. 
In this case, either $A_Z$ is non-split, or $A$ is split, in which case it does admit improper similitudes. Therefore, the additional condition~\eqref{addcond.eq} is always satisfied when the involutions have trivial discriminant. This explains why it does not occur in previous papers such as~\cite{Tignol:hyderabad} and~\cite{QT-Arason}. 
\end{enumerate}
\end{remark}

\begin{example}
\label{orthogonalhermitianform}
Let $(D,\theta)$ be a division algebra with orthogonal involution, and consider two hermitian forms $h$ and $h_0$ on a finite dimensional $D$-module $M$ with values in $(D,\theta)$. If the hermitian form $h\perp(-h_0)$ has trivial discriminant and trivial Clifford invariant, we may consider the relative Arason invariant \[e_3(h/h_0)=e_3\bigr(h\perp(-h_0)\bigl)\in H^3(F)/F^\times\cdot[D],\] as defined by Bayer and Parimala in~\cite{BP}, see also~\cite[Def. 2.11]{QT-Arason}. We claim that if the relative Arason invariant of the corresponding adjoint involutions is well defined, then we have 
\[e_3^{\ad_{h_0}}(\ad_h)=e_3(h/h_0) \in N_{\alpha,\beta}^3(F).\] In other words, $e_3^{\ad_{h_0}}(\ad_h)$ is the image of $e_3(h/h_0)$ under the natural projection 
\[H^3(F)/F^\times\cdot [D] \rightarrow H^3(F)/\bigl(F^\times\cdot [D]+\cores_{Z_0/F}(Z_0^\times\cdot[\C(\End_D(M),\ad_{h_0})])\bigl).\]
This is an easy consequence of the definition of both invariants, see notably~\cite[Rem. 2.12]{QT-Arason}. 
\end{example}

\begin{remark}
\label{sim.rem}
Reversing the viewpoint, the previous example shed new light on condition~\eqref{addcond.eq}. 
Indeed, consider two hermitian forms $h$ and $h_0$ on $M$ with values in $(D,\theta)$. Let $A=\End_D(M)$ and assume that the Clifford algebras of $\ad_h$ and $\ad_{h_0}$ are $F$-isomorphic. 
%Assume in addition that $e_3^{\ad_{h_0}}(\ad_h)$ is well defined, that is either $D_{Z}$ is non split, or $(A,\ad_{h})$ admits improper similitudes. 
By~\cite[Lemma 2]{LT}, there exists $\lambda\in F^\times$ such that  $\qform{\lambda}h\perp(-h_0)$ has trivial Clifford invariant, so that $e_3\bigl(\qform{\lambda}h\perp(-h_0)\bigr)\in H^3(F)/F^\times\cdot[D]$ is well defined. To define a relative invariant $e_3^{\ad_{h_0}}(\ad_h)$ from this element, one has to check it does not depend on the choice of the scalar $\lambda$. 
 Hence, consider $\mu\in F^\times$ such that $ \qform{\mu}h\perp(-h_0)$ has trivial Clifford invariant. 
%By the previous example, under condition~\eqref{addcond.eq}, the invariants $e_3\bigl(\qform{\lambda}h\perp(-h_0)\bigr)$ and $e_3\bigl(\qform{\mu}h\perp(-h_0)\bigr)$ in $H^3(F)/F^\times\cdot[D]$ differ by an element of \[\cores_{Z_0/F}(Z_0^\times\cdot[\C(A,\ad_{h_0})])/F^\times\cdot[D]=\cores_{Z/F}(Z^\times\cdot[\C(A,\ad_{h})])/F^\times\cdot[D],\] 
%where the two groups coincide since $\ad_h$ and $\ad_{h_0}$ have $F$-isomorphic Clifford algebras. 
%This observation can be directly obtained as follows. 
The difference between $\qform{\lambda}h\perp (-h_0)$ and $\qform{\mu}h\perp (-h_0)$ also has trivial Clifford invariant, and is Witt equivalent to $\qform{\lambda}\qform{1,-\lambda\mu}h$. 
Therefore, by Tao's formula for the Clifford algebra of a tensor product, we have either $(\lambda\mu,\disc(h))=0$ or $(\lambda\mu,\disc(h))\sim D$. In the first case, there exists $y\in Z^\times$ such that $\lambda\mu=N_{Z/F}(y)$. 
Therefore, by~\cite[Prop. 2.6]{QT-Arason}, we get 
\[e_3(\qform{\lambda}h\perp(-h_0))=e_3(\qform{\mu}h\perp(-h_0))+\cores_{Z/F}\bigr(y\cdot\C(A,\ad_{h})\bigl)\in H^3(F)/F^\times\cdot[D].\]
Assume now that $(\lambda\mu,\disc(h))\sim D$. In particular, $D_{Z}$ is split. Hence, assuming~\eqref{addcond.eq} holds, we get that the algebra with involution $(A,\ad_{h})$ admits improper similitudes. Pick such a similitude and denote by $\alpha$ its multiplier. We have $h\simeq \qform{\alpha}h$ and $D\sim (\alpha,\disc(h))$ see~\cite[(12.20)\&(13.38)]{KMRT}. Therefore, $(\lambda\mu\alpha,\disc(h))=0$, and there exists $y\in Z^\times$ such that $\lambda\mu=\alpha N_{Z/F}(y)$. Hence, we have 
\[\qform{1,-\lambda\mu}h\simeq \qform{1,-\alpha N_{Z/F}(y)}h\simeq \qform{1,-N_{Z/F}(y)}h,\] and applying again~\cite[Prop. 2.6]{QT-Arason}, we get 
\[e_3(\qform{\lambda}h\perp(-h_0))=e_3(\qform{\mu}h\perp(-h_0))+\cores_{Z/F}\bigr(y\cdot\C(A,\ad_{h})\bigl)\in H^3(F)/F^\times\cdot[D].\]
Hence in both cases, under condition~\eqref{addcond.eq}, the image of $e_3(\qform{\lambda}h\perp(-h_0))$ in $N^3_{\alpha,\beta}(F)$ does not depend on the choice of $\lambda$ such that $\qform{\lambda}h\perp(-h_0)$ has trivial Clifford invariant. This provides as an alternate definition of the relative Arason invariant $e_3^{\ad_{h_0}}(\ad_h)$. 
\end{remark}

The definition of the relative $f_3$ invariant given in~\cite[Def. 2.15]{QT-Arason} also extends to this broader setting. 
We proceed as in~\cref{f3.section}. From the fundamental relations given in~\cite[(9.12)]{KMRT}, one may easily check that for all \[y\in F^\times\cdot[A]+\cores_{Z_0/F}\bigl(Z_0^\times\cdot[\C(A,\sigma_0)]\bigr),\] we have 
\[
2 y=0\mbox{ if $\ell$ is even, and }
2y\in N_{Z_0/F}(Z_0^\times)\cdot[A]\mbox{ if $\ell$ is odd.}\\
\]
Therefore, for any algebra $A$ endowed with two orthogonal involutions $\sigma_0$ and $\sigma$ for which $e_3^{\sigma_0}(\sigma)$ is well-defined, the invariant 
\[f_3^{\sigma_0}(\sigma)=2 c\in\left\{\begin{array}{ll}
H^3(F)&\mbox{ if $\ell$ is even,}\\
H^3(F)/N_{Z_0/F}(Z_0^\times)\cdot[A]&\mbox{ if $\ell$ is odd,}\\
\end{array}\right.\]
also is well defined, where $c\in H^3(F)$ is any cohomology class such that \[e_3^{\sigma_0}(\sigma)=c\mod F^\times\cdot[A]+\cores_{Z_0/Z}\bigl(Z_0^\times\cdot[\C(A,\sigma_0)]\bigr).\]

\subsection{Orthogonal quadratic extensions of unitary involutions}
\label{orthogonalquadraticextension.sec}

Orthogonal quadratic extensions of a unitary algebra with involution are defined as symplectic ones, see~\cref{symplecticquadraticextension.sec}. 
Namely, the algebra with orthogonal involution $(A,\sigma)$ is called a quadratic extension of the algebra with unitary involution $(B,\tau)$ if $A$ contains a $\sigma$-stable quadratic \'etale $F$-algebra $\tilde F$ such that 
\[(\tilde F,\sigma_{\mid \tilde F})\simeq (F',\iota)\mbox{ and }\bigr(Z_A(\tilde F),\sigma_{\mid Z_A(\tilde F)}\bigl)\simeq (B,\tau),\] where $Z_A(\tilde F)$ denotes the centraliser of $\tilde F$ in $A$. When this holds, $A$ is a central simple $F$-algebra of degree $2n$, where $n$ is the degree of $B$, and it satisfies $A_{F'}\sim B$. In particular, this implies $B$ has exponent at most $2$. 

As in the symplectic case, the involution $\sigma$ is uniquely determined by its restriction $\tau$. More precisely, given a degree $n$ central simple algebra with $F'/F$-unitary involution $(B,\tau)$ such that $B$ has exponent at most $2$, and a central simple algebra $A$ of degree $2n$ over $F$ such that $A_{F'}\sim B$, there exists a unique orthogonal involution $\sigma$ of $A$ such that $(A,\sigma)$ is a quadratic extension of $(B,\tau)$. This is an easy consequence of the Skolem-Noether theorem, see~\cite[Prop. 1.9]{GaQ:deg12}. We call $\sigma$ the orthogonal extension of $\tau$ to $A$. 

\begin{example}
\label{splitquadext.ex}
Assume $F'=F(\sqrt\delta)$ is a field and $B$ is split, that is $B\simeq\End_{F'}(V)$ for some $n$-dimensional vector space $V$ over $F'$. The involution $\tau\simeq\ad_h$ is adjoint to a hermitian form $h:\,V\times V\rightarrow (F',\iota)$, and we denote by $q_h$ the corresponding trace, defined on the $F$ vector-space $V$ by $q_h(x)=h(x,x)$. The underlying bilinear form satisfies \[b_{q_h}(x,y)=\frac12\bigl(h(x,y)+h(y,x)\bigr),\ \mbox{for all }x,y\in V.\]
From this, one may easily check that $\ad_{q_h}$ is the orthogonal extension of $\ad_h$ to $A=\End_F(V)$. 
Indeed, multiplication by $\sqrt \delta$ is an element of $\End_F(V)$ with centraliser $B=\End_K(V)$, and it generates a subfield $\tilde F\subset \End_F(V)$ which satisfies the required conditions. 

\end{example}

The discriminant and the Clifford algebra of the orthogonal quadratic extension $(A,\sigma)$ of an algebra with unitary involution $(B,\tau)$ may be computed as follows :
\begin{lem}
\label{disc.lem}
Let $(A,\sigma)$ be a quadratic orthogonal extension of the algebra with $F'/F$-unitary involution $(B,\tau)$. 
Recall $F'=F[X]/(X^2-\delta)$ and denote by $n$ the degree of $B$, so that $A$ has degree $2n$. The invariants of the involution $\sigma$ are given by 
\begin{multline*}
d(\sigma)=\left\{\begin{array}{ll} 
\delta\in F^\times/F^{\times 2}&\mbox{ if $n=2m+1$ is odd,}\\
1\in F^\times/F^{\times 2}&\mbox{ if $n=2m$ is even.}\\
\end{array}\right.\\
\mbox{and }  
\C(A,\sigma)=\left\{\begin{array}{ll} 
0\in\br(F')&\mbox{ if $n=2m+1$ is odd,}\\
\C_+\times\C_-\in\br(F)\times\br(F)&\mbox{ if $n=2m$ is even,}\\
\end{array}\right.
\end{multline*}
with either $\C_+$  or $\C_-$ Brauer-equivalent to the discriminant algebra $\cd(\tau)$.
\end{lem} 
\begin{proof} 
If $F'\simeq F\times F$, that is $\delta=1\in F^\times/F^{\times 2}$, then both $\tau$ and $\sigma$ are hyperbolic, see~\cite[Example 1.12]{GaQ:deg12}. Therefore all invariants are trivial in this case, and the lemma holds. Hence we  may assume $F'=F(\sqrt\delta)$ is a field. 
Assume in addition that the algebra $A$ is split, so that $B\sim A_{F'}$ also is. 
The result then follows from~\cref{splitquadext.ex} together with \eqref{discJacobson} and \eqref{CliffordJacobson}. 
In the general case, consider the function field $F_A$ of the Severi-Brauer variety of $A$. The result holds after extending scalars to $F_A$, since $A$ is split over this field. Moreover, $F$ is quadratically closed in $F_A$, and by Amitsur's theorem, the kernel of the map $\br(F)\rightarrow\br(F_A)$ is $\{0,[A]\}$. In view of the fundamental relations~\cite[(9.12)]{KMRT}, this proves the lemma when $n=2m$ is even.
Assume now $n$ is odd. Since $B\sim A_{F'}$, $B$ has exponent $2$ and odd degree, hence it is split. 
Again by Amitsur's theorem, $\br(F')\rightarrow \br(F'\otimes_F F_A)$ has trivial kernel, and this concludes the proof. 
% If $m$ is odd, then $A_{F'}\sim B$ is split, since $B$ has exponent $2$. Therefore, again by Amitsur's theorem, $\br(F')\rightarrow \br(F'\otimes_F F_A)$ has trivial kernel, and this concludes the proof. 
\end{proof}
With this in hand, we now prove the following :
\begin{prop}\label{Unit-orth-1}
Let $(B,\tau_0,\tau)$ be a degree $n$ algebra with two $F'/F$ unitary involutions. We assume $B$ has exponent dividing $2$, and if $n$ is even, we assume $\tau_0$ and $\tau$ have isomorphic discriminant algebras. Let $A$ be a degree $2n$ algebra over $F$ such that $A_{F'}\sim B$, and denote by $\sigma_0$, respectively $\sigma$, the orthogonal extension of $\tau_0$, respectively $\tau$, to the algebra $A$. 
If $n$ is odd, we assume in addition that $(A,\sigma)$ admits improper similitudes. 
Then both relative Arason invariants are defined and we have 
\[e_3^{\sigma_0}(\sigma)=e_3^{\tau_0}(\tau)\mod F^\times\cdot [A].\]
\end{prop} 
\begin{remark}
\label{orthogonalextension.rem}
(1) By~\cref{disc.lem}, the involutions $\sigma_0$ and $\sigma$ have discriminant $1$ if $n$ is even and $\delta$ if $n$ is odd. It follows that condition~(\ref{addcond.eq}) is always satisfied if $n$ is even, see~\cref{relorth.rem} (3). If $n$ is odd, then $B\sim A_{F'}$ is split, since it has exponent $2$ and odd degree. Therefore, condition~(\ref{addcond.eq}) is satisfied in this case if and only if $(A,\sigma)$ admits improper similitudes, see~\cref{relorth.rem} (1).

%By~\cref{disc.lem}, the involutions $\sigma_0$ and $\sigma$ have discriminant $1$ if $n$ is even and $\delta$ if $n$ is odd. In addition, if $n$ is odd, then $B\sim A_{F'}$ is split, since it has exponent $2$ and odd degree. Therefore, condition~(\ref{addcond.eq}) is always satisfied if $n$ is even; if $n$ is odd, it is satisfied if and only if $(A,\sigma)$ admits improper similitudes, see~\cref{relorth.rem} (1) and (3).

(2) If $n$ is odd, then $B$ is split and $e_3^{\tau_0}$ has values in $N^3_{0,0}(F)=H^3(F)$. Moreover, since $\C(A,\sigma)=0\in\br(F')$, $e_3^{\sigma_0}$ has values in $H^3(F)/F^\times\cdot [A]$. The result in this case means that $e_3^{\tau_0}(\tau)$ has image $e_3^{\sigma_0}(\sigma)$ under the natural map 
\[H^3(F)\rightarrow H^3(F)/F^\times\cdot [A].\]
Assume now that $n$ is even, so that $A$ has degree divisible by $4$. Since $A_{F'}\sim B$, by the projection formula, we have $\cores_{F'/F}(F'^\times \cdot[B])=N_{F'/F}(F'^\times)\cdot[A]$. 
Therefore, the invariant $e_3^{\tau_0}$ has values in 
\[H^3(F)/\bigl(N_{F'/F}(F'^\times)\cdot[A]+F^\times\cdot[\cd(\tau_0)]\bigr).\]
On the other hand, combining~\cref{disc.lem} and~\cite[(9.12)]{KMRT}, we have 
\[\cores_{Z_0/F}(Z_0^\times\cdot[\C(A,\sigma_0)])=F^\times\cdot[\C_+]+F^\times\cdot[\C_-]=F^\times\cdot[\cd(\tau_0)]+F^\times\cdot[A].\]
Therefore, the invariant $e_3^{\sigma_0}$ has values in 
\[H^3(F)/\bigl(F^\times\cdot[A]+F^\times\cdot[\cd(\tau_0)]\bigr).\]
The result in this case means that 
$e_3^{\tau_0}(\tau)$ has image $e_3^{\sigma_0}(\sigma)$ under the natural map 
\[H^3(F)/\bigr(N_{F'/F}(F'^\times)\cdot[A]+F^\times\cdot[\cd(\tau_0)]\bigl)\rightarrow H^3(F)/\bigr(F^\times\cdot[A]+F^\times\cdot[\cd(\tau_0)]\bigl).\]
\end{remark}

\begin{proof}[Proof of Proposition~\ref{Unit-orth-1}.]
The proof is similar to the proof of~\cref{symplecticextension}. 
Since $(A,\sigma_0)$ is a quadratic orthogonal extension of $(B,\tau_0)$, there is a canonical embedding \[i:\,SU(B,\tau_0)\rightarrow \SO(A,\sigma_0),\] which factors through a homomorphism $i':\,SU(B,\tau_0)\rightarrow \Spin(A,\sigma_0)$. 
Pick a cohomology class $\xi\in H^1\bigl(F,\SU(B,\tau_0)\bigr)$ with image in $H^1\bigl(F,\PGU(B,\tau_0)\bigr)$ corresponding to $[B,\tau,\id_{F'}]$. 
Such a class exists since $\tau_0$ and $\tau$ have isomorphic discriminant algebras if $n$ is even. 
By~\eqref{ssym.eq}, $\xi$ corresponds to the class of a pair $(s,z)\in \SSym(B,\tau_0)$, where $s\in\Sym(B,\tau_0)^\times$ satisfies $\tau=\Int(s)\circ\tau_0$. Clearly, $\Sym(B,\tau_0)\subset\Sym(A,\sigma_0)$. 
Therefore, $\Int(s)\circ\sigma_0$ is an orthogonal involution on $A$, acting as $\tau$ on $B$. By uniqueness of the orthogonal extension of $\tau$ to $A$, we thus have $\sigma=\Int(s)\circ \sigma_0$. 
Hence, it follows from the definition of the relative Arason invariants in the unitary and orthogonal cases that
\[e_3^{\tau_0}(\tau)=\rho_{\SU(B,\tau_0)}(\xi)\mbox{ and }e_3^{\sigma_0}(\sigma)=\rho_{\Spin(A,\sigma_0)}\bigl({i'}^{(1)}(\xi)\bigr),\] 
where each equality holds in the relevant quotient of $H^3(F)$, see~\cref{orthogonalextension.rem}(2). 
The following lemma finishes the proof : 
\begin{lem}
The homomorphism $i':\,SU(B,\tau_0)\rightarrow \Spin(A,\sigma_0)$ has Rost multiplier $n_{i'}=1$, that is for all $\xi\in H^1\bigl(F,\SU(B,\tau_0)\bigr)$, we have 
\[\rho_{\Spin(A,\sigma_0)}\bigl({i'}^{(1)}(\xi)\bigr)=\rho_{\SU(B,\tau_0)}(\xi).\]
\end{lem}
By~\cite[Prop. 7.9(4) and Example 7.15]{GMS}, the composition 
\[\Spin(A,\sigma_0)\rightarrow \SO(A,\sigma_0)\rightarrow \SL_1(A).\]
has Rost multiplier $2$. 
Therefore, composing further on the left with $i'$, we get a map 
\[j:\,\SU(B,\tau_0)\rightarrow \SL_1(A)\] with Rost multiplier $n_{j}=2 n_{i'}$ by~\cite[Prop.  7.9(1)]{GMS}. 
On the other hand, we may also compute the Rost multiplier of $j$ directly, assuming both groups are split. Under this assumption, the involutions are hyperbolic, and we may identify $(A,\sigma_0)$ with 
$M_{2n}(F)$ endowed with the involution defined by  
\[\sigma_0\begin{pmatrix}A&B\\C&D\end{pmatrix}=\begin{pmatrix}D^t&B^t\\C^t&A^t\end{pmatrix},\mbox{ for all }A,B,C,D\in M_n(F).\]
Therefore, as in the symplectic case, the map 
\[\psi:\,\bigl(M_n(F)\times M_n(F)^{\mathrm{op}},\varepsilon)\rightarrow (M_{2n}(F),\sigma_0),\ \ 
(X,Y^{\mathrm{op}})\mapsto \begin{pmatrix}X&0\\0&Y^t\end{pmatrix}\]
preserves the involutions. Moreover, $\tilde F=F[\Delta]$, where $\Delta=\begin{pmatrix}1&0\\0&-1\end{pmatrix}$, satisfies 
$(\tilde F,{\sigma_0}_{\mid\tilde F})\simeq (F\times F, \iota)$, and the centraliser, endowed with the restriction of $\sigma_0$, coincides with the image of $\psi$. Hence, $\sigma_0$ is the orthogonal extension of $\varepsilon$ to $M_{2n}(F)$. Identifying $\SU\bigl(M_n(F)\times M_n(F)^{\mathrm{op}},\varepsilon\bigr)$ with $\SL_n(F)$ by $X\mapsto (X,(X^{-1})^{\mathrm {op}})$, we get that $j$ corresponds in the split case to the map 
\[\,\SL_n(F)\rightarrow \SL_{2n}(F), \ \ X\mapsto \begin{pmatrix} X&0\\0&(X^{-1})^{t}\end{pmatrix}.\]
By~\cref{rostmultiplier}, we have $n_j=2n_{i'}=2$ and this concludes the proof. 
\end{proof}

\begin{remark}
If $n$ is odd, Proposition~\ref{Unit-orth-1} has a much easier proof. Under this assumption, 
the algebra $B$ is split, and $A$ is similar to a quaternion algebra $Q$. 
By~\cite[Prop. 2.2]{Peyre}, $H^3(F)/F^\times\cdot[Q]$ maps injectively to $H^3(F_Q)$, where $F_Q$ is the function field of the Severi-Brauer variety of $Q$. Therefore, since $e_3^{\sigma_0}$ has values in $H^3(F)/F^\times\cdot[Q]$ and $F_Q$ is a generic splitting field for $Q$, it is enough to check the result when $Q$ is split, in which case it is a direct consequence of~\cref{splitquadext.ex} and~\cref{split}. 
\end{remark}

\subsection{Unitary involutions with an orthogonal descent}
\label{orthogonaldescent.sec}
In this section, we consider a central simple algebra $B$ with two $F'/F$-unitary involutions $\tau_0$ and $\tau$, and we assume they both admit an orthogonal descent to the same algebra $A$ over $F$, that is 
\begin{equation}
(B,\tau_0)\simeq (A,\sigma_0)\otimes_F (F',\iota)\mbox{ and }(B,\tau)\simeq (A,\sigma)\otimes_F (F',\iota),\
\end{equation}
for some orthogonal involutions $\sigma$ and $\sigma_0$ of $A$. 
As in the symplectic case, this does not hold in general, see~\cref{existencedescent}. 
When it is satisfied, we aim at comparing the invariants of the pairs $(\tau_0,\tau)$ and $(\sigma_0,\sigma)$ of unitary and orthogonal involutions. 

We first prove the following: 
\begin{prop}
\label{splitorthogonaldescent.prop}
Let $\tau_0$ and $\tau$ be two unitary involutions of the degree $n$ split algebra $B=\End_{F'}(V)\simeq \End_F(V)\otimes_F F'$. If $n$ is even, we assume in addition that the discriminant algebras $\cd(\tau_0)$ and $\cd(\tau)$ are isomorphic. 
Then, there exists two quadratic forms $q_0$ and $q$ over $F$, having the same discriminant $d\in F^\times/F^{\times 2}$, and such that 
$\tau_0\simeq \ad_{q_0}\otimes\iota$ and $\tau\simeq \ad_q\otimes\iota$. 
Moreover, the relative Arason invariant of the unitary involutions $\tau_0$ and $\tau$ is given by 
\[e_3^{\tau_0}(\tau)=(\delta)\cdot[\C(q-q_0)]\in\left\{\begin{array}{ll}
H^3(F)&\mbox{ if $n$ is odd},\\
H^3(F)/F^\times\cdot(\delta, d)&\mbox{ if $n$ is even},\\
\end{array}\right.\]
where $[C(q-q_0)]\in \br(F)$ denotes the Brauer class of the full Clifford algebra of $q-q_0$. 
\end{prop} 

\begin{proof}
Let $h_0$ be a hermitian form defined on $V$ and with values in $(F',\iota)$ such that $\tau_0\simeq\ad_{h_0}$. 
Pick a diagonalisation $h_0\simeq\qform{a_1,\dots, a_n}$. Then the elements $a_i\in F'$ are fixed under $\iota$, hence they belong to $F$. So, one may consider the $n$ dimensional quadratic form $q_0\simeq \qform{a_1,\dots, a_n}$, and a direct computation shows $\tau_0\simeq\ad_{q_0}\otimes\iota$. Define similarly a quadratic forms $q$ such that $\tau\simeq \ad_q\otimes\iota$. 
If $n$ is odd, replacing $q$ by $\qform{d(q_0)d(q)^{-1}}q$, we may assume that $q$ and $q_0$ have the same discriminant, denoted by $d$. 
Assume now that $n$ is even. Since $\cd(\tau_0)$ and $\cd(\tau)$ are isomorphic,~(\ref{DiscOrthDescent.eq}) shows there exists $z\in F'^\times$ such that $d(q_0)=d(q)N_{F'/F}(z)$. The one dimensional hermitian forms $\qform{a_1}$ 
and $\qform{\iota(z)a_1z}\simeq\qform{N_{F'/F}(z)a_1}$ are isomorphic, so we may replace $q_0$ by $\qform{a_1N_{F'/F}(z),a_2,\dots, a_n}$ and we get $d(q)=d(q_0)\in F^\times/F^{\times 2}$ as required. Hence in both cases, we have $d(q_0)=d(q)$, so that $q-q_0$ has trivial discriminant. 
We may now compute $e_3^{\tau_0}(\tau)$ using~\cref{split}. Since the trace forms of $h_0$ and $h$ respectively are $q_{h_0}=\pform{\delta}\otimes q_0$ and $q_h=\pform{\delta}\otimes q$, we get 
\[e_3^{\tau_0}(\tau)=e_3\bigl(\pform{\delta}\otimes(q-q_0)\bigr)=(\delta)\cdot e_2(q-q_0)=(\delta)\cdot [\C(q-q_0)],\] 
and this concludes the proof. 
\end{proof}

\begin{remark}
If $n$ is odd, and if $\tau_0$ admits an orthogonal descent, then the algebra $B$ has exponent dividing $2$, hence it is split. 
This case is covered by the previous proposition. Therefore, we will assume in the following that $n$ is even. 
\end{remark}

If the algebra $B$ is non split, we assume a stronger condition on the orthogonal descent, namely that their relative Arason invariant is well defined; we get : 
\begin{prop}
\label{orthdescent.prop}
Let $n=2m$ be an even integer, and let $B$ be a degree $n$ algebra endowed with two unitary $F'/F$ involutions $\tau_0$ and $\tau$. We assume $\tau_0$ and $\tau$ admit orthogonal descents to the same algebra, that is  
\[B\simeq A\otimes _F F',\ \tau_0\simeq \sigma_0\otimes \iota\mbox{ and }\tau\simeq \sigma\otimes\iota,\]
for some central simple algebra $A$ over $F$, and some orthogonal involutions $\sigma_0$ and $\sigma$ of $A$. 
We assume in addition that the Clifford algebras $\C(A,\sigma_0)$ and $\C(A,\sigma)$ are $F$-isomorphic, and $(A,\sigma)$ satisfies condition~(\ref{addcond.eq}), so that $e_3^{\sigma_0}(\sigma)$ and $f_3^{\sigma_0}(\sigma)$ are well defined. 
Then, the involutions $\tau_0$ and $\tau$ have isomorphic discriminant algebras, and their relative Arason invariant is given by  
\[e_3^{\tau_0}(\tau)=f_3^{\sigma_0}(\sigma)\in H^3(F)/N_{F'/F}(F'^\times)\cdot[A]+F^\times\cdot[\cd(\tau_0)].\]
\end{prop}

\begin{remark} 
(1) In the situation of the proposition above, the computation of the
discriminant algebra recalled in~(\ref{DiscOrthDescent.eq}) shows that $e_3^{\tau_0}$ has values in 
\[\left\{\begin{array}{ll} 
H^3(F)/N_{F'/F}(F'^\times)\cdot[A]+F^\times\cdot\bigr(\delta,d(\sigma_0)\bigl)&\mbox{ if $m$ is even, and }\\
H^3(F)/F^\times\cdot[A]+F^\times\cdot\bigr(\delta,d(\sigma_0)\bigl)\cdot&\mbox{ if $m$ is odd.}\\
\end{array}\right.\]
On the other hand, $f_3^{\sigma_0}(\sigma)$ belongs to the group 
$H^3(F)$ if $m$ is even, and the quotient $H^3(F)/N_{Z_0/F}(Z_0^\times)\cdot[A]$ if $m$ is odd, see the end of \cref{orthogonalinvariant.sec}. 
The statement of the proposition means that $e_3^{\tau_0}(\tau)$ is the image of $f_3^{\sigma_0}(\sigma)$ under the relevant projection. 

(2) Note that Propositions~\ref{splitorthogonaldescent.prop} and ~\ref{orthdescent.prop} are compatible when they both  apply. Indeed, assume $n=2m$ even, $B$ is split, $\tau_0\simeq\ad_{q_0}\otimes\iota$ and $\tau\simeq \ad_q\otimes\iota$, where $q_0$ and $q$ have the same discriminant $d$ and isomorphic even Clifford algebras $\C_0(q_0)$ and $\C_0(q)$. 
By the structure theorem for Clifford algebras~\cite[Thm. 2.5]{Lam}, $\C_0(q_0)$ is the centralizer in $\C(q_0)$ of the quadratic extension $Z=F[X]/(X^2-d)$. Therefore, by~\cite[Chap. 8, Thm. 5.4]{Scharlau}, $\C({q_0}_Z)$ is Brauer equivalent to $\C_0(q_0)$. The same argument also applies to $q$. 
Hence, $\C({q_0}_Z)$ and $\C(q_Z)$ are isomorphic. 
So there exists $a\in F^\times$ such that $[\C(q-q_0)]=(d,a)$, and we get $(\delta)\cdot[\C(q-q_0)]=(\delta,d,a)\in H^3(F)$, and $e_3^{\tau_0}(\tau)=0\mod F^\times\cdot(\delta, d)$. On the other hand, since the Arason invariant for quadratic forms takes values in $_2H^3(F)$, we also have $f_3^{\sigma_0}(\sigma)=0$ in this situation. 
\end{remark} 

\begin{proof}[Proof of Proposition~\ref{orthdescent.prop}.] 
By assumption, the algebras $\C(A,\sigma_0)$ and $\C(A,\sigma)$ are $F$-isomorphic, so they have isomorphic centers. Therefore, $\sigma_0$ and $\sigma$ have the same discriminant, and it follows that $\cd(\tau_0)$ and $\cd(\tau)$ are isomorphic. Hence $e_3^{\tau_0}(\tau)$ is well defined. 

Moreover, since $(B,\tau_0)\simeq (A,\sigma_0)\otimes_F (F',\iota)$, there is a natural inclusion \[\SO(A,\sigma_0)\subset \SU(B,\tau_0),\] which induces a map $j:\,\Spin(A,\sigma_0)\rightarrow \SU(B,\tau_0)$. When both groups are split, identifying $\SU(M_n(F)\times M_n(F)^{\mathrm{op}},\varepsilon)$ with $\SL_n(F)$, this map corresponds to the canonical composition 
\[\Spin_n\rightarrow \SOn\rightarrow \SL_n.\]
Since $n=2m$ is even, it has Rost multiplier $2$ by~\cite[Example 7.15]{GMS}, hence for all $\xi\in H^1\bigr(F,\Spin(A,\sigma_0)\bigl)$, we have 
\[\rho_{\SU(B,\tau_0)}\bigr(j^{(1)}(\xi)\bigl)=2\rho_{\Spin(A,\sigma_0)}(\xi)\in H^3(F).\]
Let $\xi\in H^1\bigr(F,\Spin(A,\sigma_0)\bigl)$ be a cohomology class with image in $H^1(F,\PGO^+\bigl(A,\sigma_0)\bigr)$ corresponding to the triple $(A,\sigma,\psi)$, so that 
\[e_3^{\sigma_0}(\sigma)=\rho_{\Spin(A,\sigma_0)}(\xi)\in H^3(F)/\bigl(F^\times\cdot[A]+\cores_{Z_0/F}(Z_0^\times\cdot[\C(A,\sigma_0)])\bigl),\]
and 
\[f_3^{\sigma_0}(\sigma)=2\rho_{\Spin(A,\sigma_0)}(\xi)\in \left\{\begin{array}{ll}
H^3(F)&\mbox{ if $m$ is even,}\\
H^3(F)/N_{Z_0/F}(Z_0^\times)\cdot [A]&\mbox{ if $m$ is odd.}\\
\end{array}\right.\]
The image of $\xi$ in $H^1\bigl(F,\SO(A,\sigma_0)\bigr)$ corresponds under the canonical bijection described in~\cite[(29.27)]{KMRT} to the class of a pair $(s,z)$, where $s\in \Sym(A,\sigma_0)$ satisfies $\sigma=\Int(s)\circ \sigma_0$. We thus have $\tau=\Int(s\otimes 1)\circ\tau_0$, and it follows that $j^{(1)}(\xi)\in H^1\bigl(F,\SU(B,\tau_0)\bigr)$ has image in $H^1\bigr(F,\PGU(B,\tau_0)\bigl)$ corresponding to $(B,\tau,\id_{F'})$. Therefore, 
$e_3^{\tau_0}(\tau)$ coincides with $\rho_{\SU(B,\tau_0)}(j^{(1)}(\xi))$ in the relevant quotient of $H^3(F)$, and this concludes the proof. 
\end{proof}

\section{Specific results in small degree} 
\label{lowdegree.section}

In this section, we explore the relative Arason invariant for unitary involutions on small degree algebras. In particular, we prove it is classifying in degree $4$ under some additional assumption. In degree $4$ and $6$, the hyperbolic Arason invariant characterizes hyperbolic involutions. This in not true anymore if the underlying algebra has degree $8$ and index at most $4$. In this situation, we define a new absolute Arason invariant, taking a totally decomposable involution as a base point instead of the hyperbolic one, and we compare it with the hyperbolic Arason invariant. It is not known whether or not this invariant detects totally decomposable involutions, but partial results in this direction are established. 

\subsection{Unitary involutions on degree $4$ algebras}

Throughout this section, we let $B$ be a degree $4$ central simple algebra over $F'$ endowed with two $F'/F$-unitary involutions $\tau_0$ and $\tau$. We assume in addition that $\tau_0$ and $\tau$ have isomorphic discriminant algebra, denoted by $A$. Hence, the relative Arason invariant $e_3^{\tau_0}(\tau)$ is well defined, and it belongs to $N^3_{\alpha,\beta}(F)$, where $\alpha=[A]\in \br(F)$ and $\beta=[B]\in \br(F')$. 
The main result in this subsection is the following : 
\begin{thm}
\label{fund-degree4}
Let $(B,\tau_0,\tau)$ be a degree $4$ algebra with two $F'/F$-unitary involutions with isomorphic discriminant algebras. 
\begin{enumerate}
\item[(a)] If $e_3^{\tau_0}(\tau)=0\in N^3_{\alpha,\beta}(F)$
then $(B,\tau)$ and $(B,\tau_0)$ are $F$-isomorphic, that is $(B,\tau)$ is $F'$-isomorphic either to $(B,\tau_0)$ or to its conjugate $(^\iota B,^\iota \tau_0)$. 
\item[(b)] Assume in addition that either $B$ and $^\iota B$ are not isomorphic, or $(B,\tau_0)$ and $(^\iota B, ^\iota \tau_0)$ are $F'$-isomorphic. Then  $e_3^{\tau_0}(\tau)=0\in N^3_{\alpha,\beta}(F)$ if and only if $\tau$ and $\tau_0$ are isomorphic. 
\end{enumerate}
\end{thm}

\begin{remark}
If $B$ is split and has degree $4$, Merkurjev introduced a degree $3$ invariant denoted by $e_3(B,\tau)$ in~\cite[\S 4]{Me2000}. It can be interpreted as a relative invariant as follows. Since $B$ is split, $\tau=\ad_h$ for some rank $4$ hermitian form $h$ with discriminant $d(h)\in F^\times/N_{F'/F}(F'^\times)$. Consider the isotropic hermitian form $h_0=\qform{1,-d(h)}+{\mathbb {H}}$.  Both involutions $\tau_0$ and $\tau$ have discriminant algebra Brauer equivalent to the quaternion algebra $\bigl(\delta,d(h)\bigr)$ and we have $e_3(B,\tau)=e_3^{\ad_{h_0}}(\ad_h)$. In particular, it follows from our theorem that this invariant vanishes if and only if $\tau$ is isotropic, which is Lemma 4.1 in Merkurjev's paper. 
\end{remark}

\begin{proof}[Proof of Theorem~\ref{fund-degree4}.] 
Let $(B,\tau_0,\tau)$ be as in the statement, and assume the relative Arason invariant
$e_3^{\tau_0}(\tau)=0\in N^3_{\alpha,\beta}(F).$ 
As explained in \S\ref{relative.section}, $e_3^{\tau_0}(\tau)$ is the class in $N^3_{\alpha,\beta}(F)$ of $\rho_{\SU(B,\tau_0)}(\xi)$, where $\xi\in H^1\bigl(F,\SU(B,\tau_0)\bigr)$ has image $\eta\in H^1\bigl(F,\PGU(B,\tau_0)\bigr)$ corresponding to the class of the triple $[B,\tau,\id_{F'}]$. Hence, under our assumption, there exists $\mu\in F'^\times$ such that 
\begin{equation}
\label{trivial}
\rho_{\SU(B,\tau_0)}(\xi)=\cores_{F'/F}\bigl((\mu)\cdot[B]\bigr)\mod F^\times\cdot [A].
\end{equation}
In order to prove assertion (a), we use the equivalence of groupoids ${\sf A}_3\equiv {\sf D}_3$ established in~\cite[(15.D)]{KMRT}, and we interpret the element of $H^3(F)/F^\times\cdot[A]$ above as the Arason invariant of some hermitian form of orthogonal type. 

Specifically, denote by $\sigma_0$ the canonical involution of the discriminant algebra $A$ of $\tau_0$. Thus, $A$ is a degree $6$ central simple algebra over $F$, $\sigma_0$ is an orthogonal involution of $A$, and by~\cite[(15.D)]{KMRT}, the Clifford algebra of $(A,\sigma_0)$ is $F$-isomorphic to $(B,\tau_0)$. Similarly, the canonical involution of the discriminant algebra of $\tau$ is an orthogonal involution $\sigma$ of $A$, and the Clifford algebra of $(A,\sigma)$ is $F$-isomorphic to $(B,\tau)$. The algebra $A$ is Brauer equivalent to a division algebra $D$, which is either $F$ or an $F$-quaternion algebra. Pick an orthogonal involution $\theta$ of $D$ and a hermitian module $(M,h_0)$ over $(D,\theta)$ such that $(A,\sigma_0)\simeq (\End_D(M),\ad_{h_0})$. 

By~\cite[(15.26)]{KMRT}, we have $\PGU(B,\tau_0)\simeq \PGO^+(A,\sigma_0)\simeq \PGO^+(h_0)$, hence we may view $\xi$ and $\eta$ as elements of $H^1\bigl(F,\Spin(h_0)\bigr)$ and $H^1\bigl(F,\PGO^+(h_0)\bigr)$, respectively. The morphism $\Spin(h_0)\rightarrow \PGO^+(h_0)$ factors through $\SO(h_0)$. Let $\zeta\in H^1\bigr(F,\SO(h_0)\bigl)$ be the image of $\xi$ under the induced map. The cocycle $\zeta$ corresponds to a hermitian form $h$ over $M$, with image $\eta$ in $H^1\bigr(F,\PGO^+(h_0)\bigl)$, so that $\ad_h\simeq \sigma$. Moreover, since $\zeta$ is the image of $\xi\in H^1\bigr(F,\Spin(h_0)\bigl)$, we may consider the relative Arason invariant 
\[e_3(h/h_0)=\rho_{\Spin(A,\sigma_0)}(\xi)\mod F^\times\cdot[A],\]

which coincides with $e_3(h\perp (-h_0))$, see \cref{orthogonalhermitianform}. Hence~\eqref{trivial} gives  
\[e_3(h\perp (-h_0))=\cores_{F'/F}\bigl((\mu)\cdot[B]\bigr)\mod F^\times\cdot [A].\]
On the other hand, the Clifford algebra of $(\ad_{h_0})_{F'}$ is $B\otimes_F F'\simeq B\times B$. 
By~\cite[Prop. 2.6]{QT-Arason}, it follows that 
\[e_3(\qform{1,-N_{F'/F}(\mu)}h_0)=\cores_{F'/F}\bigl((\mu)\cdot[B]\bigr)\mod F^\times\cdot [A].\]
Hence, adding these two equalities, and using~\cite[Lemma 2.2]{QT-Arason}, we get 
\[e_3\bigr(h\perp\qform{-N_{F'/F}(\mu)}h_0\bigl)=2\cores_{F'/F}\bigl((\mu)\cdot[B]\bigr)\mod F^\times\cdot [A].\]
By the fundamental relations given in~\cite[(9.12)]{KMRT}, we have $2[B]=[A_{F'}]$, hence 
\[2\cores_{F'/F}\bigl((\mu)\cdot[B]\bigr)=\cores_{F'/F}\bigl((\mu)\cdot[A_{F'}]\bigr)=N_{F'/F}(\mu)\cdot [A],\] and this proves $e_3\bigr(h\perp\qform{-N_{F'/F}(\mu)}h_0\bigl)=0\in H^3(F)/F^\times\cdot [A]$. 
Hence the corresponding adjoint involution of $\End_D(M\oplus M)$ also has trivial Arason invariant. The algebra $\End_D(M\oplus M)$ has degree $12$ and index at most $2$. So we may apply~\cite[Thm. 4.1]{QT-Arason}, and we get that the involution $\ad_{h\perp\qform{-N_{F'/F}(\mu)}h_0}$ is hyperbolic, and so is the hermitian form. 
Consequently, the forms $h_0$ and $h$ are similar, so $\sigma_0$ and $\sigma$ are isomorphic, and their Clifford algebras $(B,\tau_0)$ and $(B,\tau)$ are $F$-isomorphic. This proves (a), and assertion (b) follows immediately. 
\end{proof}

\begin{remark}
(1) As observed by Garibaldi-Petersson in~\cite{GP}, the additional condition in (b) is satisfied if and only if the group $\SU(B,\tau)$ admits an outer automorphism defined over $F$ as soon as there is no Tits class obstruction, see also~\cite[Prop. 2.1]{QT-AutExt}. 
By \cref{relorth.rem}(2), this holds if and only if if $(A,\sigma_0)$ satisfies~\eqref{addcond.eq}. 
Specifically, since $A$ has degree $6$, the canonical homomorphism 
\[c:\,\PGO(A,\sigma_0)\rightarrow {\mathrm {Aut}}_F(B,\tau_0)\] is an isomorphism, see~\cite[(15.1)]{KMRT}. Moreover, by~\cite[(13.2)]{KMRT}, improper similitudes of $(A,\sigma_0)$ correspond to automorphisms of $(B,\tau_0)$ that act as $\iota$ on $F'$. Hence, $(A,\sigma_0)$ admits improper similitudes if and only if $(B,\tau_0)$ and $(^\iota B, ^\iota \tau_0)$ are $F'$-isomorphic. 

(2) When \eqref{addcond.eq} is satisfied, the Arason invariant $e_3^{\sigma_0}(\sigma)$ is well defined. Moreover, identifying the groups $\PGU(B,\tau_0)$ and $\PGO^+(A,\sigma_0)$, one may check that the cocycle $\eta$ corresponding to the class of $[B,\tau,\id_{F'}]$ also corresponds to $[A,\sigma,\varphi]$ for some isomorphism $\varphi$ between the center of the Clifford algebras of $(A,\sigma_0)$ and $(A,\sigma)$, which are both isomorphic to $F'$. It follows that 
\[e_3^{\tau_0}(\tau)=e_3^{\sigma_0}(\sigma)\in N^3_{\alpha,\beta}(F).\]

Therefore, the argument in the proof of \cref{fund-degree4} also shows that if $A$ is a degree $6$ algebra with two orthogonal involutions $\sigma_0$ and $\sigma$ with $F$-isomorphic Clifford algebras and such that  \eqref{addcond.eq} holds, then $e_3^{\sigma_0}(\sigma)=0\in N^3_{\alpha,\beta}(F)$ if and only if $\sigma$ and $\sigma_0$ are isomorphic. 
\end{remark} 

In particular, the theorem shows that $e_3^{\hyp}$ detects hyperbolic involutions : 
\begin{cor}
Assume $B$ has degree $4$ and index at most $2$, and $\tau$ is an $F'/F$-unitary involution. Then $e_3^{\hyp}(\tau)$ vanishes if and only if $\tau$ is hyperbolic. 
\end{cor}
\begin{proof}
Let $\tau_0$ be the hyperbolic involution of $B$; then $^\iota \tau_0$ also is hyperbolic. Hence the additional condition in~\cref{fund-degree4}(b) is automatically satisfied and this proves the corollary. 
\end{proof}

\subsection{Unitary involutions on degree $8$ and index $4$ algebras}

The main tool in this section is the cohomological invariant introduced in~\cite{Bar14}. 
Let $D$ be a biquaternion algebra over $F'$ with trivial $F'/F$-corestriction. The algebra $D$ always has a descent up to Brauer equivalence, but it generally does not have a descent up to isomorphism, see~\cite[Rem. 4.8]{Bar14} and~\cref{hypnotdec}. Specifically, the Brauer class of $D$ always comes from a Brauer class of $F$, but the underlying division algebra might be of degree $8$. In~\cite[Def. 4.2, Prop. 4.4]{Bar14}, a cohomological invariant \[\delta_{F'/F}(D)\in H^3(F,\mu_2)/\cores_{F'/F}(F'^\times\cdot [D])\] is associated to $D$, and it is proved that $\delta_{F'/F}(D)=0$ if and only if $D=D_0\otimes_F F'$ for some biquaternion algebra $D_0$ defined over $F$. In other words, the invariant $\delta_{F'/F}(D)$ detects whether or not $D$ has a descent to $F$ up to isomorphism. Since $H^3(F,\mu_2)$ embeds canonically into $H^3(F)$, we may as well consider $\delta_{F'/F}(D)$ as an element of $H^3(F)/\cores_{F'/F}(F'^\times\cdot [D])=N_{0,\beta}(F)$, with $\beta=[D]\in\br(F')$. 

Throughout this section, $(B,\tau)$ denotes a central simple algebra over $F'$ of degree $8$ and index at most $4$, endowed with an $F'/F$-unitary involution. Hence, $B\simeq M_2(D)$ for some biquaternion algebra $D$ over $F'$. Moreover, $B$ and $D$ have trivial corestriction, so that $\delta_{F'/F}(D)$ is well defined. We first prove the following : 

\begin{prop}
\label{totdecdeg8}
Let $(B,\tau)$ be an algebra with $F'/F$-unitary involution. We assume $B=M_2(D)$ for some biquaternion algebra $D$. If $(B,\tau)$ is totally decomposable, then 
\[\cd(\tau)\sim 0\mbox{ and }e_3^{\hyp}(\tau)=\delta_{F'/F}(D)\in H^3(F)/\cores_{F'/F}(F'^\times\cdot [B]).\]
\end{prop} 

\begin{proof}
Assume $(B,\tau)$ satisfies the conditions of the Proposition. In particular, it is totally decomposable, that is there exist quaternion algebras with unitary involutions $(Q_i,\tau_i)$ over $F'$ such that 
\[(B,\tau)\simeq \otimes_{i=1}^3(Q_i,\tau_i).\] By~\cite[(2.22)]{KMRT}, each $(Q_i,\tau_i)$ decomposes as $(H_i, \gamma_i)\otimes (F',\iota)$, where $H_i$ is a quaternion algebra over $F$, and $\gamma_i$ its canonical involution. Hence $(B,\tau)$ admits the algebra with symplectic involution \[(C,\gamma)\simeq (H_1,\gamma_1)\otimes (H_2,\gamma_2)\otimes(H_3,\gamma_3)\] as a symplectic descent. In particular, its discriminant algebra is split by \cref{symplecticdescent.prop}. Moreover, we may use the same proposition to compute $e_3^{\hyp}(\tau)$. Indeed, we have : 
\begin{lem} 
The hyperbolic unitary involution $\tau_0$ of $B$ admits a symplectic descent to $C$. 
\end{lem} 
This can be checked as follows. If $C$ is not division, it has a symplectic hyperbolic involution $\gamma_0$ and by uniqueness of the hyperbolic involution of a given type on a given algebra, we have $\tau_0\simeq \gamma_0\otimes \iota$. 
Assume now $C$ is division. Since $C\otimes_F F'\simeq B\simeq M_2(D)$ is not division, $C$ contains a quadratic subfield $\tilde F$ which is isomorphic to $F'$, and the centralizer of $\tilde F$ in $C$ is isomorphic to $D$. In particular, it has degree $4$. Applying~\cite[(4.14)]{KMRT} to the subalgebra $(\tilde F, \id_{\tilde F})$ of $(C,\gamma)$, we get a symplectic involution $\gamma_0$ of $C$ which acts as $\id_{\tilde F}$ on $\tilde F$. Therefore there exists $u\in C$ such that $u^2=\delta$ and $\gamma_0(u)=u$. Consider the $F'/F$-unitary involution $\gamma_0\otimes \iota$ on $B$ and the idempotent $e=\frac{1}{2}(1+\delta^{-2}u\otimes \sqrt{\delta})\in B$. A direct computation shows that $(\gamma_0\otimes \iota)(e)=1-e$. Hence $\gamma_0\otimes \iota$ is hyperbolic, and so $\tau_0\simeq \gamma_0\otimes \iota$.

With this in hand, Proposition~\ref{symplecticdescent.prop} applies to $(B,\tau_0,\tau)$ and we get : 
\[e_3^{\hyp}(\tau)=e_3^{\tau_0}(\tau)=e_3^{\gamma_0}(\gamma)=e_3(\gamma_0)-e_3(\gamma)\in H^3(F)/\cores_{F'/F}(F'^\times\cdot{[B]}).\]
Moreover, since $(C,\gamma)$ is totally decomposable, $e_3(\gamma)=0$ by~\cite[Thm. B]{GPT}, and it remains to prove that $e_3(\gamma_0)=\delta_{F'/F}(D)\in H^3(F)/\cores_{F'/F}(F'^\times\cdot{[B]}).$ If $C$ is not division, then $D$ has a descent up to isomorphism and $\gamma_0$ is hyperbolic so $\delta(D)=0$ and $e_3(\gamma_0)=0$. If $C$ is division, the definition of $\gamma_0$ above shows it contains a subfield $\tilde F$ isomorphic to $F'$ over which $\gamma_0$ acts as identity. The result follows in this case from~\cite[Rem. 4.3]{Bar14}, and this finishes the proof.  
\end{proof} 

\begin{example} 
When $D$ has a descent up to isomorphism, all totally decomposable involutions have trivial $e_3^{\hyp}$ invariant, even the non hyperbolic ones. Specifically, 
for any quaternion algebras $H_1$ and $H_2$ over $F$, we have \[\delta_{F'/F}\bigl((H_1\otimes_F H_2)\otimes_F F'\bigr)=0.\] Hence, given $F$-linear involutions $\rho_i$ of $H_i$ for $i=1,2$, and $\lambda\in F^\times$, we get that 
\[(B,\tau)=(M_2(F),\ad_{\pf{\lambda}})\otimes_F (H_1,\rho_1)\otimes_F (H_2,\rho_2)\otimes_F (F',\iota)\] 
satisfies $e_3^{\hyp}(\tau)=0$. Nevertheless, $\tau$ is not hyperbolic in general, so the invariant $e_3^{\hyp}$ does not characterize hyperbolic involutions in degree $8$. To construct an explicit example, let $F$ be a purely transcendental extension $F=k(t)$ of a field $k$, $H_1\otimes H_2$ a biquaternion algebra over $k$ which is division over $k'=k(\sqrt{\delta})$, and pick $\lambda=t$. The involution $\tau$ is adjoint to the hermitian form $\qform{1,-t}$ with values in $(H_1\otimes H_2\otimes F',\rho_1\otimes \rho_2\otimes \iota)$, which is anisotropic. Note that the algebra with involution considered here is a 'generic sum' of two copies of $\gamma_1\otimes\gamma_2\otimes \iota$ in the sense of~\cite{QT-AutExt}. 
\end{example} 

As a consequence of the previous proposition, we have : 
\begin{cor}
\label{hypversustotdec}
Let $B=M_2(D)$ for some biquaternion algebra $D$ over $F'$ with trivial $F'/F$-corestriction. The hyperbolic unitary involution $\tau_0$ is totally decomposable if and only if the invariant $\delta_{F'/F}(D)$ vanishes. 
\end{cor}

\begin{proof}
Assume first that $\tau_0$ is totally decomposable. By the proposition above, we have $e_3^{\hyp}(\tau_0)=\delta_{F'/F}(D)$. On the other hand, since $\tau_0$ is hyperbolic,  $e_3^{\hyp}(\tau_0)=0$, and this proves the first implication. 

To prove the converse, assume $\delta_{F'/F}(D)=0 \in H^3(F)/ \cores_{F'/F}(F'^\times\cdot [D])$. By~\cite[Prop. 4.4]{Bar14}, there exists a biquaternion algebra $D_0$ over $F$ such that $D\simeq D_0\otimes_F F'$. 
By uniqueness of the hyperbolic unitary involution of $B$, we have
\[(B,\tau_0)\simeq \bigr((M_2(F),\gamma)\otimes_F (D_0,\rho)\bigl)\otimes _F (F',\iota),\] where $\gamma$ is the canonical involution of $M_2(F)$, which is hyperbolic, and $\rho$ is any involution of the first kind on $D_0$. Since $D_0$ is a biquaternion algebra, it admits decomposable involutions, and the result follows immediately. 
\end{proof}

\begin{example}
\label{hypnotdec}
Let $C\simeq H_1\otimes H_2\otimes H_3$ be a totally decomposable algebra over $F$ such that $F'\subset C$ and there is no quaternion subalgebra of $C$ containing $F'$. Explicit examples are given in~\cite[Cor. 4.5]{Bar16}. Since $F'$ is a subfield of $C$, $C_{F'}$ is isomorphic to $M_2(D)$ for some biquaternion algebra $D$ over $F'$. Assume $D=D_0\otimes_F F'$ for some biquaternion algebra $D_0$ over $F$. Then $C\simeq D_0\otimes H$ for some quaternion algebra $H$ containing $F'$. Hence, by construction of $C$, such a $D_0$ does not exist and it follows that $\delta_{F'/F}(D)\not =0$. 
As a consequence, the hyperbolic unitary involution on $M_2(D)$ is not totally decomposable and the totally decomposable involution $\tau=\gamma_1\otimes\gamma_2\otimes\gamma_3\otimes\iota$ of $C_{F'}\simeq M_2(D)$ satisfies $e_3^{\hyp}(\tau)=\delta_{F'/F}(D)\not=0$. 
\end{example}
In view of this example, in order to characterize totally decomposable involutions using a cohomological invariant, it seems relevant to work with a totally decomposable base point rather than the hyperbolic one. 
The following lemma explains how this can be done : 
\begin{lem}
\label{e3dec}
Let $\tau_1$ be a totally decomposable unitary involution of the degree $8$ algebra $M_2(D)$. The discriminant algebra of $\tau_1$ is split, and for all unitary involutions $\tau$ of $B$ with split discriminant algebra, we have : 
\[e_3^{\tau_1}(\tau)=e_3^{\hyp}(\tau)+\delta_{F'/F}(D)\in H^3(F)/\cores_{F'/F}(F'^\times\cdot[B]).\]
In particular, $e_3^{\tau_1}(\tau)$ does not depend on the choice of a totally decomposable involution $\tau_1$, and will be denoted by $e_3^{\mathrm{td}}(\tau)$. 
\end{lem}
\begin{proof}
The discriminant algebra of $\tau_1$ is split by~\cref{totdecdeg8}. 
Let $\tau_0$ be the hyperbolic unitary involution of $B$. By \cref{Chasles} and~\cref{order}, we have 
\[e_3^{\tau_1}(\tau)=e_3^{\tau_1}(\tau_0)+e_3^{\tau_0}(\tau)=e_3^{\hyp}(\tau_1)+e_3^{\hyp}(\tau).\]  
This finishes the proof since $e_3^{\hyp}(\tau_1)=\delta_{F'/F}(D)$ by \cref{totdecdeg8}.  
\end{proof} 

Hence, we have defined another degree $3$ invariant for unitary involutions with split discriminant algebra by 
\[e_3^{\mathrm{td}}(\tau)=e_3^{\tau_1}(\tau)\in H^3(F)/\cores_{F'/F}(F'^\times\cdot[B]),\] where $\tau_1$ is an arbitrary totally decomposable involution of $B$. 
This invariant is related to the hyperbolic Arason invariant by 
\[e_3^{\mathrm{td}}(\tau)=e_3^{\hyp}(\tau)+\delta_{F'/F}(D),\] and the two coincide if and only if $\delta_{F'/F}(D)=0$. 
Moreover, since any totally decomposable involution can be chosen as a base point, we have : 
\begin{lem}
\label{dec.st}
Let $\tau$ be a totally decomposable unitary involution of $M_2(D)$, where $D$ is a biquaternion algebra over $F'$. Then $\tau$ has split discriminant algebra and 
\[e_3^{\mathrm{td}}(\tau)=0.\]
\end{lem} 
It is not known whether or not the converse holds in general. In the remaining part of this section, we prove partial results in this direction. The first one deals with the split case : 
\begin{prop}\label{dec-split-case}
Let $\tau$ be an $F'/F$-unitary involution on a split degree $8$ algebra $B=M_8(F')$.  The involution $\tau$ is totally decomposable if and only if \[\cd(\tau)\mbox{ is split and }e_3^{\mathrm{td}}(\tau)=0 \in H^3(F).\]
\end{prop}
\begin{remark}
If $B$ is split, so is $D$. Hence, $\delta_{F'/F}(D)=0$ and $e_3^{\mathrm{td}}=e_3^{\hyp}$. Therefore we could replace $e_3^{\mathrm{td}}$ by $e_3^{\hyp}$ in this case. Moreover, both invariants have values in $H^3(F)/\cores(F'^\times\cdot[B])=H^3(F)$. 
\end{remark} 
\begin{proof}[Proof of Proposition~\ref{dec-split-case}.]
One implication is given by \cref{dec.st}, so we only need to prove the converse. 
Assume $e_3^{\hyp}(\tau)=0 \in H^3(F)$. Let $h$ be a hermitian forms over $(F',\iota)$ such that $\tau=\ad_{h}$. As observed in \cref{splithyp.example}, we have $e_3^{\hyp}(\tau)=e_3(q_h)$. Since $q_h$ is a $16$-dimensional quadratic form, the equality $e_3(q_h)=0$ implies that $q_h$ is similar to a $4$-fold Pfister form $\varphi$ over $F$. Moreover, since $q_h$ contains the subform $\pf{\delta}$, the form $\varphi_{F'}$ is isotropic, hence hyperbolic, since it is a Pfister form. It follows that the pure subform of $\varphi$ represents $-\delta$, see ~\cite[X.  Thm. 4.5]{Lam}. Hence, by  ~\cite[X. Thm. 1.5]{Lam},  $\varphi$ is of the form $\pf{\delta, a, b,c}$ for some $a,b,c\in F^\times$. Therefore $h$ is similar to $\pf{a, b,c}$, now viewed as a hermitian form over $(F',\iota)$ and $\tau\simeq \ad_\pf{a, b,c}$ which is totally decomposable.
\end{proof}
Another situation where involutions with trivial $e_3^{\mathrm{td}}$ are totally decomposable is the following. Using orthogonal decompositions of the underlying hermitian forms, one may always view a unitary involution $\tau$ of $M_2(D)$ as an orthogonal sum of two involutions $\rho_1$ and $\rho_2$ of $D$. In other words, $\tau=\theta_1$, given by 
\begin{equation}
\label{sumbis.eq}
\theta_1\begin{pmatrix} x&y\\z&t\\\end{pmatrix}=\begin{pmatrix} 
\rho_1(x)&-\rho_1(z)s^{-1}\\
-s\rho_1(y)&s\rho_1(t)s^{-1}\\
\end{pmatrix},\mbox{ for all }x,y,z,t\in D,
\end{equation}
for some well chosen $s\in \Sym(D,\rho_1)^\times$ such that $\rho_2=\Int(s)\circ\rho_1$. This also means $\tau$ is the the adjoint involution with respect to the hermitian form $\qform{1,-s}$ with values in $(D,\rho_1)$. This hermitian form represents $1$; this can only be achieved for certain choices of $\rho_1$. 
The next result deals with unitary involutions for which we may find such a decomposition with $\rho_1$ having split discriminant algebra. 
\begin{prop}
Let $B=M_2(D)$, where $D$ is a biquaternion algebra over $F'$ and assume $D$ is endowed with a unitary involution $\rho_1$ with split discriminant algebra. For all $s\in  \Sym(D,\rho_1)^\times$, the involution $\theta_1$ of $M_2(D)$ given by \eqref{sumbis.eq} is totally decomposable if and only if \[\cd(\theta_1)\mbox{ is split and }e_3^{\mathrm{td}}(\theta_1)=0\in H^3(F)/\cores_{F'/F}\bigl(F'^\times\cdot[B]\bigr).\]
\end{prop}
\begin{proof}
One implication is given by~\cref{totdecdeg8}; let us prove the converse. 
Since $D$ has degree $4$ and $\rho_1$ has split discriminant algebra, there exists quaternion algebras $H_1$ and $H_2$ over $F$ with canonical involutions $\gamma_1$ and $\gamma_2$, respectively, such that \[(D,\rho_1)\simeq (H_1,\gamma_1)\otimes_F (H_2,\gamma_2)\otimes_F (F',\iota),\]
see~\cite[thm. 3.1]{KQM} and~\cite[(2.22)]{KMRT}. 
In particular, $D$ has a descent up to isomorphism, so $\delta_{F'/F}(D)=0$ and $e_3^{\mathrm{td}}(\theta_1)=e_3^{\hyp}(\theta_1)$. Moreover, $\id_{H_1\otimes H_2}\otimes \iota$ is an $F$ automorphism of $(D,\rho_1)$ acting as $\iota$ on $F'$, so that $(D,\rho_1)\simeq (^\iota D,^\iota\rho_1)$.

By assumption, $\theta_1$ has split discriminant algebra. Hence, the same computation as in the proof of \cref{sum} shows $\cd(\theta_1)\sim \bigl(\delta,\Nrd_{D}(s)\bigr)\sim\cd(\rho_2)$. Therefore $\rho_2$ also has split discriminant algebra. Hence, \cref{sum} applies to $(D,\rho_1,\rho_2)$ and we get that $e_3^{\mathrm{td}}(\theta_1)=e_3^{\hyp}(\theta_1)=e_3^{\rho_1}(\rho_2)=0$. By \cref{fund-degree4}(b), we get that $\rho_2$ is isomorphic to $\rho_1$. So, there exists $u\in D^\times$ such that $\rho_2= \Int(u)\circ\rho_1 \circ \Int(u^{-1})=\Int(u\rho_1(u))\circ\rho_1$. It follows that $s=\lambda u\rho_1(u)$ for some $\lambda\in F^\times$. So the hermitian form $\qform{1,-s}$ with values in $(D,\rho_1)$ is isomorphic to $\qform{1,-\lambda}$ and a direct computation shows 
\[
(M_2(D),\theta_1)\simeq (M_2(F),\ad_{\pf{\lambda}})\otimes_F (H_1,\gamma_1)\otimes_F (H_2,\gamma_2)\otimes_F (F',\iota).
\] 
In particular, it is totally decomposable and this finishes the proof. 
\end{proof}

\subsection{Unitary involutions on degree $6$ algebras}
In this section, we assume that the algebra $B$ has degree $6$ and is endowed with a hyperbolic $F'/F$-unitary involution, so that $B=M_2(D)$ where $D$ has degree $3$.
We first prove that $e_3^{\hyp}$ detects hyperbolicity in this case : 
\begin{prop}\label{property-degree6}
\label{hyp.deg6}
Let $B$ be a degree $6$ and index dividing $3$ central simple $F$-algebra with center $F'$ and let $\tau$ be an $F'/F$-unitary involution on $B$. 
The involution $\tau$ is hyperbolic if and only if $\mathcal D(\tau)$ is split and $e_3^{\hyp}(\tau)=0\in N^3_{0,\beta}(F)$.
\end{prop}
\begin{proof}
One implication follows from~\ref{hyperbolic.section}, so we only need to prove the converse. 
Assume $\cd(\tau)$ is split and $e_3^{\hyp}(\tau)=0\in N^3_{0,\beta}(F)$, and let us prove $\tau$ is hyperbolic. As we did in degree $8$, using a diagonalisation of the underlying hermitian form, we may view $\tau$ as an orthogonal sum of two unitary involutions $\rho_0$ and $\rho$ of $D$. By~\cref{sum}, we have $e_3^{\hyp}(\tau)=e_3^{\rho_0}(\rho)=0$. Therefore, since $D$ has degree $3$, $\rho_0\simeq \rho$ by~\cref{deg3.2}. 
Hence we may assume $s=1$ in~\eqref{sum.eq}, so that $\tau$ is adjoint to the hermitian form $\qform{1,-\lambda}$ with values in $(D,\rho_0)$, where $\lambda\in F^\times$. 
Moreover, by the same computation as in the proof of~\cref{sum}, the discriminant algebra of $\tau$ is $(\delta,\lambda^3)=0$. So, there exists $x\in F'^\times$ such that $\lambda=N_{F'/F}(x)=x\iota(x)=x\rho_0(x)$, and this proves that the hermitian form $\qform{1,-\lambda}$ over $(D,\rho_0)$ is isomorphic to $\qform{1,-1}$. We get that the involution $\tau$ is hyperbolic as expected. 
\end{proof}

The relative Arason invariant also is classifying under some additional conditions. Let us first assume the underlying algebra $B$ is split. 
\begin{prop}\label{degree6-split}
Let $B=M_6(F')$ and consider two $F'/F$-unitary involutions $\tau_0$ and $\tau$ of $B$. We assume that $\mathcal D(\tau_0)=\mathcal D(\tau)=0$.  The involutions $\tau_0$ and $\tau$ are isomorphic if and only if $e_3^{\tau_0}(\tau)=0 \in H^3(F)$.
\end{prop}
\begin{proof} Assume $e_3^{\tau_0}(\tau)=0\in H^3(F)$. 
Since $B$ is split, there are $6$-dimensional hermitian forms $h_0$ and $h$ with values in $(F', \iota)$ and with trivial discriminant such that $\tau_0=\ad_{h_0}$ and $\tau=\ad_h$. So, we may find $p,q,r,s,t\in F$ and $u\in F'^\times$ such that
\[
h =\qform{p,q,r,s,t, pqrst N_{F'/F}(u)} = \qform{pqr}\qform{pq,pr,qr, pqrs, pqrt,st N_{F'/F}(u)}.
\] 
The $3$-dimensional quadratic forms $\qform{pq, pr, qr}$ and $\qform{pqrs, pqrt, stN_{F'/F}(u)}$ have respective determinant $1$ and $N_{F'/F}(u)$. Hence, $h$ is similar to the form 
\begin{eqnarray*}
\qform{a,b, ab} -\qform{N_{F'/F}(u)}\qform{c,d,cd}
\end{eqnarray*}
for some  $a,b,c,d\in F^\times$. 
Hence, the Jacobson trace of $h$ is
\[
q_h\simeq \pf{\delta}\left(\qform{a,b,ab} -\qform{N_{F'/F}(u)}\qform{c,d,cd}\right)\simeq \pf{\delta}\qform{a,b,ab,-c,-d,-cd}, 
\]
which is the product of the $1$-fold Pfister form $\pf{\delta}$ with an Albert form. 
Likewise, we may find $a_0, b_0, c_0, d_0\in F^\times$ such that
\[
q_{h_0}\simeq \pf{\delta}\qform{a_0,b_0,a_0b_0,-c_0,-d_0,-c_0d_0}. 
\]
Since $e_3^{\tau_0}(\tau)=e_3(q_h- q_{h_0})\in H^3(F)$ by Lemma \ref{split}, the assumption implies that $e_3(q_h- q_{h_0})=0$. Hence, $q_h= q_{h_0} \bmod I^4 (F)$. By~\cite{Hoff} this implies that the two forms $q_h$ and $q_{h_0}$ are similar. Consequently, the corresponding hermitian forms $h$ and $h_0$ also are similar, and the involutions $\tau_0\simeq \ad_{h_0}$ and $\tau\simeq \ad_h$ are isomorphic. 
\end{proof}
%%%
The next result extends slightly \cref{hyp.deg6}. More precisely, we assume $B=M_2(D)$ for some degree $3$ algebra $D$, and we consider two involutions having orthogonal sum decompositions with a common factor. In other words, both are adjoint to hermitian forms representing $1$ with values in $(D,\rho_0)$ for some $F'/F$-unitary involution $\rho_0$ of $D$. 
\begin{prop}
Let $(D,\rho_0)$ be a degree $3$ central simple $F$-algebra with an $F'/F$-unitary involution, and let $s_1, s_2\in \Sym(D,\rho_0)^\times$. Consider the unitary involutions $\tau_1$ and $\tau_2$ respectively adjoint to the hermitian forms ${\pf{s_1}}$  and ${\pf{s_2}}$ with values in $(D,\rho_0)$. We assume in addition that $\mathcal D(\tau_1)$ and 
$\cd(\tau_2)$ are both split. Then $\tau_1$ and $\tau_2$ are isomorphic if and only if $e_3^{\tau_1}(\tau_2)=0\in N_{0,\beta}^3(F)$. 
\end{prop}
\begin{proof}
The argument is pretty similar to the proof of~\cref{hyp.deg6}; the main difference is that we need to assume here that the two involutions admit orthogonal sum decompositions with a common factor $\rho_0$, while this is always the case when one of them is hyperbolic. More precisely, consider $\tau_1$ and $\tau_2$ as in the statement, and assume $e_3^{\tau_1}(\tau_2)=0$. Consider the unitary involutions of $D$ defined by $\rho_i=\Int(s_i)\circ \rho$ for $i=1,2$. 
By~\cref{Chasles}, \cref{order} and ~\cref{sum}, we have the following equalities in $H^3(F)/N_{F'/F}(F'^\times\cdot [B])=H^3(F)/N_{F'/F}(F'^\times\cdot [D])$: 
\[0=e_3^{\tau_1}(\tau_2)=e_3^{\hyp}(\tau_1)+e_3^{\hyp}(\tau_2)=e_3^{\rho_0}(\rho_1)+e_3^{\rho_0}(\rho_2)=e_3^{\rho_1}(\rho_2).\]
Therefore, by~\cref{deg3.2}, the involutions $\rho_1$ and $\rho_2$ are isomorphic. Hence, there exists $\lambda\in F^\times$ such that $s_2=\lambda s_1$. Comparing the discriminant algebras of $\tau_1$ and $\tau_2$, the same kind of computation as in the proof of~\cref{hyp.deg6} shows $\lambda=x\rho_0(x)$ for some $x$ and this concludes the proof. 
\end{proof}

\end{document}